\documentclass[11pt, reqno]{amsart}
\usepackage[utf8]{inputenc}
\usepackage{a4wide,amsfonts,amsfonts, amsmath,amssymb,amsthm,epsfig,cite,graphicx,hyperref,
color, esint,fancyhdr, enumerate, latexsym,amsrefs, makecell, float, blindtext, multicol,ulem}
\usepackage{calligra}
\usepackage{tcolorbox}
\usepackage{microtype}

\usepackage[T1]{fontenc} 
\usepackage{newtxtext,newtxmath}

\newtheorem{thm}{Theorem}[section]
\newtheorem{cor}[thm]{Corollary}
\newtheorem{lem}[thm]{Lemma}
\newtheorem{prop}[thm]{Proposition}

\numberwithin{equation}{section}

\newcommand{\mbb}{\mathbb}
\newcommand{\ra}{\rightarrow}
\newcommand{\z}{\zeta}
\newcommand{\pa}{\partial}
\newcommand{\ov}{\overline}
\newcommand{\sm}{\setminus}
\newcommand{\ep}{\epsilon}

\newcommand{\al}{\alpha}
\newcommand{\Om}{\Omega}

\newcommand{\la}{\lambda}

\newcommand{\de}{\delta}

\newcommand{\mf}{\mathbb}
\newcommand{\om}{\omega}
\newcommand{\ga}{\gamma}
\newcommand{\Ga}{\Gamma}
\newcommand{\be}{\beta}

\begin{document}

\title{A Quasiconformal Squeezing function and Fridman invariant and applications}

\author{Diganta Borah, Steven G. Krantz, Prachi Mahajan, Kaushal Verma}

\address{DB: Indian Institute of Science Education and Research Pune, Pune 411008, India}
\email{dborah@iiserpune.ac.in}

\address{SGK: Department of Mathematics, Campus Box 1146, Washington University in St. Louis, St.
Louis, Missouri 63130, USA}
\email{sk@math.wustl.edu }

\address{PM: Department of Mathematics, Indian Institute of Technology Bombay, Powai, Mumbai 400076, India}
\email{prachi.mjn@iitb.ac.in}

\address{KV: Department of Mathematics, Indian Institute of Science, Bengaluru-560012, India}
    \email{kverma@iisc.ac.in}

\begin{abstract}
    The squeezing function and Fridman invariant are two examples of biholomorphic invariants that have proved to be useful in exploring questions in complex geometry. Motivated by their construction, we define conformal analogues of these invariants by using suitable families of $K$-quasiconformal mappings, and study several basic properties such as continuity, domain dependence and boundary behaviour near spherically extreme points. It is also shown that these conformal invariants are explicitly related. As an application, we propose and prove variants of the Ferrand--Obata theorem in conformal geometry that apply to domains in $\mbb R^n$. This highlights the close analogy that exists between the Ferrand--Obata and Wong--Rosay theorems.
\end{abstract}

\maketitle

\section{Introduction}

\noindent The principle that abundant symmetry leads to rigidity appears repeatedly in geometry. Notable examples include the Myers-Steenrod theorem, relating the geometry of a manifold to its isometry group; Obata's characterization of the sphere through eigenvalue conditions; Mostow's rigidity theorem, which establishes the topological determination of hyperbolic geometry in many settings; the Wong-Rosay theorem, which singles out the unit ball in $\mathbb C^n$ among strongly pseudoconvex domains admitting a non-compact holomorphic automorphism group; and finally, the Ferrand--Obata theorem which is an analogue of the Wong--Rosay theorem in the conformal category.

\medskip
 
The Ferrand--Obata theorem \cite{F1, F2, Ferrand-coll, Obata1, Obata2} provides the principal motivation for the present work. We therefore begin by recalling its setting and highlighting the principal features that are relevant to our discussion.

\medskip

Let $M$ be a smooth $n$-dimensional Riemannian manifold with metric $g$, which is not assumed to be complete. Denote by $[g]$, the conformal class of $g$. Thus, $[g]$ consists of all Riemannian metrics on $M$ of the form $\lambda g$ for some smooth positive function $\la$ on $M$. A conformal diffeomorphism $f: (M, [g]) \ra (N, [h])$ is a diffeomorphism between $M$ and $N$ such that $f^{\ast} h \in [g]$. Let ${\rm Conf}(M, g)$ be the group of all conformal self-diffeomorphisms $f: (M, [g]) \ra (M, [g])$. The Ferrand--Obata theorem asserts that if ${\rm Conf}(M, g)$ acts non-properly on $M$, then $M$ is conformally equivalent to $S^n$ or $\mbb R^n$. 

\medskip

What happens when $M = D$ is a bounded domain in $\mathbb R^n$, $n \ge 3$, equipped with $g = g_{euc}$, the standard euclidean metric? As an example, let us examine the case of $D = \mbb B^n$, the unit ball in $\mbb R^n$. By Liouville's theorem, the conformal group consists of M\"obius maps that preserve $\mbb B^n$, and the action of this group on $\mbb B^n$ is seen to be proper. The Ferrand--Obata theorem cannot be applied per se since the hypothesis on the conformal group action is not satisfied. The right hypothesis, then, for a domain-version of the Ferrand--Obata theorem seems to be non-compactness of the conformal group, which is weaker than requiring non-properness of its action. In analogy with the Wong--Rosay theorem, we are thus led to propose the following domain-version of the Ferrand--Obata theorem: {\it a smoothly bounded domain $D \subset \mbb R^n$ with non-compact ${\rm Conf}(D, g_{euc})$ is conformally equivalent to $\mbb B^n$}.
While the result may well be folklore, we have been unable to locate a suitable reference. Accordingly, we provide a brief sketch of a proof. The underlying idea is that each element of ${\rm Conf}(D, g_{euc})$ is a M\"obius map by Liouville's theorem and hence can be regarded as an element of the conformal group ${\rm Conf}(S^n, g_{sph}) \backsimeq O(n+1, 1)$; here, $g_{sph}$ denotes the spherical metric on $S^n$. Now exploit the global dynamics of the M\"obius group on $S^n$; the boundary $\pa D$ becomes an invariant compact subset of the sphere, and a divergent sequence of automorphisms that arises from the non-compactness of ${\rm Conf}(D, g_{euc})$ produces attractor-repeller dynamics on $S^n$. The rigidity comes from the fact that the only proper domains in the conformal sphere that admit such large M\"obius actions are round balls. More details are provided in $\S 7$.

\medskip

However, unlike the conformal case, holomorphic automorphisms of a given domain need not extend past the boundary or even to it in a reasonable sense. Thus, proofs of the Wong-Rosay theorem \cite{Rosay, Wong} must necessarily work from within the given domain. The scaling technique provides a paradigm to produce the desired biholomorphic equivalence with the unit ball in $\mbb C^n$, by staying within the domain, when the holomorphic automorphism group is known to be non-compact and additional assumptions on the boundary are imposed; see \cite{KV} for a $K$-quasiconformal analogue of the Wong--Rosay theorem. The interaction between non-compactness of the holomorphic automorphism group and the Levi geometry of a given domain in $\mbb C^n$ has been effectively used to show several variants (for example, \cite{BF1, BF2}) of the Wong-Rosay theorem. The survey by Isaev–Krantz \cite{IK} provides an overview of this area. The Greene--Krantz conjecture
that predicts the geometric nature of the boundary of a domain near boundary accumulation
points is a central question that remains open in its complete generality. 

\medskip

The purpose of this work is to propose an alternative way of approaching such questions, and possibly others, in the conformal category by defining a set of conformal invariants in terms of suitable $K$-quasiconformal mappings, $1 \le K < \infty$. They are modeled on the classical squeezing function (\cite{DGZ2012}) and its dual, the Fridman invariant (\cite{Fridman1, Fridman2}), which are defined in terms of a suitable class of holomorphic mappings. Both have been studied extensively and have proved to be rather effective in many questions in complex geometry. The new invariants, which we christen the $K$-squeezing function and the $K$-Fridman invariant, form an entire family parametrised by the dilatation $K \ge 1$. 

\medskip

The motivations for introducing these new classes of invariants in terms of $K$-quasiconformal mappings are two-fold. First, the numerous parallels between quasiconformal and holomorphic mappings, particularly at a large-scale or structural level, suggest the possibility of introducing new objects and viewpoints that further elucidate and strengthen the connections between these two theories. Second, Ferrand worked with the conformal capacity, which is a conformal invariant defined in terms of an energy-like functional. On the other hand, the $K$-squeezing function and $K$-Fridman invariant are defined via a variational problem that encodes a model domain and mappings from or to it. Their very definition makes them amenable to the question of uniformization in the conformal category; having a family of such invariants, instead of just one, adds further flexibility in their utility. 

\medskip

The $K$-squeezing function and $K$-Fridman invariant are defined and explored in $\S 3$ and $\S 4$ below. A comparison between them can be found in $\S 6$. The holomorphic versions of the squeezing function and Fridman invariant are known to provide alternate proofs of the Wong--Rosay theorem. In much the same spirit, their $K$-quasiconformal analogues can be used to show that a smoothly bounded domain in $\mbb R^n$ whose conformal diffeomorphism group acts transitively on it must be conformally equivalent to the ball $\mbb B^n \subset \mbb R^n$. Thus, this version of the Ferrand--Obata theorem and the Wong--Rosay theorem share not only a common conceptual framework but also notable similarities in their technical development.

\medskip

To understand the basic properties of the $K$-squeezing function and the $K$-Fridman invariant, we will need to work with Ferrand's $\mu$-pseudodistance which is related to conformal capacity (see \cite{F5}), and this is summarized in the following sections.

\section{The $\mu$-pseudodistance}
Let $M$ be a non-compact Riemannian manifold of class $C^1$ and dimension $n \geq 2$. Let $ H (M) $ denote the space of continuous real-valued
functions $ u $ on $ M $ that admit a generalized differential $  \nabla u $ with
\[
I(u,M) = \int_M | \nabla u|^n \; dv < + \infty.
\]
Here, $ dv $ denotes the Riemannian volume element on $M$, and $ | \nabla u| $ is the norm of the generalized differential with respect to the Riemannian metric on $M$. 

\medskip

Let $ H_0(M) $ be the subspace of $ H ( M ) $ consisting of functions with compact support in $ M $. For a compact set $ C \subset M $, a function $ u \in H_0(M) $ is said to be \textit{admissible for $ C $} if $ u = 1 $ on $ C $. The conformal capacity of $ C $ in $ M $ is defined by  
\[
\mbox{Cap}_M(C) = \inf_u I(u, M), 
\]
where the infimum is taken over all functions $ u \in H_0(M) $ that are admissible for $ C $. Define the function $ \mu_M: M \times M \rightarrow [0, \infty) $ by
\[
\mu_M (x,y) = \inf_{C} \mbox{Cap}_M(C),
\]
where the infimum is taken over the set of compact continua $ C \subset M $ connecting $ x $ and $ y $.

\medskip

From \cite{F5}, it is known that $ \mu_{M} $ is a pseudodistance on $M$: $\mu_M(x,x)\geq 0$, $\mu_{M}(x,x) = 0 $ for all $ x \in M $, and $\mu_M$ satisfies the triangle inequality. It was also shown in \cite{F5} that $\mu_{M}$ identically vanishes or is a distance on $M$. 

\medskip

Accordingly, Ferrand introduced the following classification of non-compact Riemannian manifolds into two complementary classes: $ M $ is said to be of Class I if $ \mu_M $ is identically zero and it is said to be of Class II if $\mu_M$ is a distance on $M$. The Euclidean space $ \mathbb{R}^n $ is an example of a Class I manifold, whereas every bounded subdomain of $ \mathbb{R}^n $ is of Class II. Moreover, every open submanifold of a Class II manifold is again in Class II. In dimension $ n=2$, Ferrand's classes coincide with the classical dichotomy for non-compact Riemann surfaces: Class I consists of parabolic Riemann surfaces, whereas Class II consists of hyperbolic Riemann surfaces.

\medskip

The pseudodistance $\mu_M$ is continuous on $ M \times M $ and it has the following quasi-invariant property: If $f: M\to M'$ is a $K$-quasiconformal homeomorphism, then
\begin{equation}\label{qinv}
K^{-1}\mu_{M}(x,y) \leq \mu_{M'}\big(f(x),f(y)\big) \leq K \mu_{M}(x,y)
\end{equation}
for all $x,y \in M$. In particular, $\mu_M$ is conformally invariant. Moreover, $\mu_{M}$ satisfies the following monotonicity property: if $M_1$ is an open submanifold of $M_2$ then
\begin{equation}\label{mon}
\mu_{M_2}(x,y) \leq \mu_{M_1}(x,y)
\end{equation}
for all $ x, y \in M_1 $. Furthermore, if $M$ is of Class II, then the topology induced by $\mu_{M}$ is the same as the manifold topology.

\medskip

From \cite[Theorem~8.6]{Vu}, we know that for the unit ball $\mathbb{B}^n$,
\begin{equation}\label{mu-ball}
\mu_{\mathbb{B}^n}(x,y)=\ga_n\left(\frac{1}{\tanh \frac{1}{2} \rho_{\mathbb{B}^n}(x,y)}\right),
\end{equation}
where $\rho_{\mathbb{B}^n}$ is the hyperbolic distance on $\mathbb{B}^n$ and is given by
\begin{equation}\label{hyp-ball}
\tanh\left(\frac{1}{2} \rho_{\mathbb{B}^n}(x,y)\right)=\frac{\vert x-y\vert}{\sqrt{\vert x-y\vert^2 +(1-\vert x \vert^2)(1-\vert y\vert^2)}},
\end{equation}
and $\ga_n(s)$, defined for $s>1$, denotes the modulus of the curve family connecting the complementary components of the Gr\"otzsch ring $R_{G,n}(s)$. Recall that $R_{G,n}(s) \subset \mathbb{R}^n$ is the ring domain whose complementary components in $\overline{\mathbb{R}^n}$ (the one-point compactification of $\mbb R^n$) are $\overline{\mathbb{B}^n}$ and the ray $[se_1,\infty]$ along the $e_1$-direction in $\mbb R^n$. The notion of the modulus of a curve family is recalled in $\S 3.1$.

\medskip

The function $\ga_n:(1, \infty) \to (0,\infty)$ is continuous, strictly decreasing, and satisfies (see, for example, \cite[Lemma~5.53]{Vu})
\begin{equation}\label{ga-dec}
\lim_{s \to 1+}\ga_n(s)=\infty \quad \text{and} \quad \lim_{s \to \infty} \ga_n(s)=0.
\end{equation}
Following \cite[$\S 7.43$]{Vu}, the function $\phi_K:(0,1) \to (0,1)$ defined by
\begin{equation}\label{qc-dist-fn-phi}
\phi_{K,n}(r)=\frac{1}{\ga_n^{-1}\left(K\gamma_n\left(\frac{1}{r}\right)\right)}
\end{equation}
is continuous, strictly increasing, and satisfies
\begin{equation}\label{phi-limit}
\lim_{r\to 0+}\phi_{K,n}(r)=0,\quad \lim_{r \to 1-}\phi_{K,n}(r)=1.
\end{equation}
We also note the following upper bound \cite[7.47]{Vu}
\begin{equation}\label{phi-ub}
\phi_{K,n}(r) \leq \la_n^{1-\al}r^{\al}, \quad \al=K^{1/(1-n)}
\end{equation}
where $\la_n$ is the Gr\"otzsch ring constant (see \cite[$\S 7.21$]{Vu} for the definition).

\medskip

\noindent When $n$ is fixed, we will denote $R_{G,n}$, $\ga_n$, and $\phi_{K,n}$ by $R_G$, $\ga$, and $\phi_K$, respectively.

\section{$K$-Squeezing function}

\noindent For suitable spaces $X, Y$ and dilatation $K \ge 1$, let $QC_K(X, Y)$ denote the family of $K$-quasiconformal embeddings from $X$ into $Y$, i.e., $f \in QC_K(X, Y)$ if $f: X \ra f(X) \subset Y$ is a $K$-quasiconformal homeomorphism. When there is no restriction on the dilatation, the family of all quasiconformal embeddings from $X$ into $Y$ will be denoted by $QC(X, Y)$. For $p \in \mathbb R^n$ and $\delta > 0$, $B(p;\delta)$ will denote the Euclidean ball around $p$ of radius $\delta$. 

\medskip

Let $\Om \subset \mf{R}^n$ be a bounded domain and fix a dilatation $K \ge 1$. The $K$-squeezing function of $\Om$ at $x \in \Om$ is
\begin{equation}\label{s^K}
s^K_{\Om}(x)=\sup\big\{r>0: B(0;r) \subset f(\Om), f \in QC_K(\Om, \mathbb B^n), f(x) = 0\big\}.
\end{equation}
Thus, $s^K_{\Om}(x)$ measures the largest Euclidean ball contained in the images of all $K$-quasiconformal embeddings from $\Om$ into the ball $\mathbb B^n$ that are normalized to map $x$ to the origin $0 \in \mathbb B^n$. Observe that $s^K_{\Omega}: \Omega \to (0,1]$.

\medskip

If $\phi:\Omega \to \Omega'$ is an $L$-quasiconformal homeomorphism, then
\[
s^K_{\Om}(x) \leq s^{KL}_{\Om'}\big(\phi(x)\big), \quad x \in \Om.
\]
Indeed, let $r$ be a candidate for $s^K_{\Omega}(x)$. Then there exists a $K$-quasiconformal map $f: \Omega \to \mathbb{B}^n$ such that $f(x)=0$ and $B(0;r) \subset f(\Omega)$. Thus, $f\circ\phi^{-1}: \Omega' \to \mathbb{B}^n$ is a $KL$-quasiconformal map satisfying $f\circ \phi^{-1}(\phi(x))=0$ and $B(0;r) \subset f\circ \phi^{-1}(\Omega')$. Therefore, $r$ is also a candidate for $s^{KL}_{\Om'}(\phi(x))$. The desired inequality follows immediately. In particular, if $L=1$, i.e., if $\phi$ is conformal, then
\[
s^K_{\Om}(x) = s^K_{\Om'}\big(\phi(x)\big), \quad x \in \Om.
\]
Thus, $s^K_{\Om}$ is a conformal invariant.


To see the dependence of $s^K_{\Om}$ on the parameter $K$, fix an $x \in \Om$, and note that the definition implies that $s^K_{\Om}$ is nondecreasing in the parameter $K$: if $K_1\leq K_2$ then
\[
s^{K_1}_{\Om}(x)\leq s^{K_2}_{\Om}(x).
\]
Therefore
\[
\lim_{K \to 1}s_{\Om}^K(x)=\inf_{K \in [1,\infty)} s_{\Om}^K(x) \quad \text{and} \quad \lim_{K \to \infty}s_{\Om}^K(x) =\sup_{K\in [1,\infty)} s_{\Om}^K(x).
\]
Later, it will be shown that $s^K_{\Omega}$ is strictly increasing in $K$ if $\Om$ is not quasiconformally equivalent to $\mbb B^n$. Furthermore, since $s^1_{\Om}(x) \leq s^K_{\Om}(x)$ for all $K \ge 1$, it follows that
\[
\inf_{K \in [1,\infty)}s_{\Om}^K(x)=\min_{K \in [1,\infty)}s_{\Om}^K(x)=s^1_{\Om}(x).
\]
On the other hand, we will show later that 
\begin{equation}\label{S}
\sup_{K\in [1,\infty)} s_{\Om}^K(x)=S_{\Om}(x):=
\sup \big\{r> 0: B(0;r) \subset f(\Om), f \in QC(\Om, \mathbb B^n), f(x) = 0 \big\}.
\end{equation}

\subsection{Existence of extremals for $s^K_{\Om}$} Since $s^K_{\Om}$ is defined via a variational problem, it is natural to ask whether there is an extremal map that realizes $s^K_{\Om}$ for a given $x \in \Om$. To see that this is indeed the case, we have to first estimate the modulus of a certain family of paths.

\medskip

\noindent Let $\overline{\mathbb{R}}^n=\mathbb{R}^n \cup \{\infty\}$ be the one-point compactification of $\mf{R}^n$. By a curve in $\overline{\mathbb{R}}^n$, we mean a continuous map $\gamma: I \to \overline{\mathbb{R}}^n$ where $I\subset \mathbb{R}$ is an interval. If $I$ is a closed interval, then the curve $\gamma$ is called a path. Let $\Gamma$ be a family of curves in $\overline{\mathbb{R}}^n$. By an \textit{admissible density} for $\Gamma$ we mean a Borel function $\rho: {\mathbb{R}}^n \to [0,\infty]$ with the property
\[
\int_{\gamma} \rho \,ds \geq 1
\]
for every locally rectifiable member $\gamma$ of $\Gamma$. We will denote the set of all admissible densities for $\Gamma$ by $\text{Adm}(\Gamma)$. The \textit{modulus} of $\Gamma$ is defined by
\[
M(\Gamma)= \inf_{\rho \in \text{Adm}(\Gamma)} \int_{\mathbb{R}^n} \rho^n \,dm
\]
with the convention that $M(\Gamma)=\infty$ if $\text{Adm}(\Gamma)=\emptyset$, which happens only if $\Ga$ contains a constant path.

\medskip

Let $\Gamma$ and $\Gamma'$ be curve families in $\overline{\mathbb{R}}^n$. We write $\Gamma < \Gamma'$, and say that $\Gamma'$ is minorized by $\Gamma$, if every $\gamma \in \Gamma'$ has a subcurve which belongs to $\Gamma$. In this case, Theorem~6.4 of \cite{V} guarantees that
\begin{equation}\label{mod-dec}
M(\Gamma) \geq M(\Gamma').
\end{equation}

If $E,F,G$ are subsets of $\overline{\mathbb{R}}^n$ and $\gamma:[a,b] \to \overline{\mathbb{R}}^n$ is a path satisfying $\gamma(a) \in E$, $\gamma(b) \in F$, and $\gamma(t) \in G$ for all $a<t<b$, then we say that $\gamma$ is a path in $G$ joining $E$ to $F$.

\begin{lem}\label{positive-capacity}
Let $\Omega \subset \mathbb{R}^n$ be a domain with
\[
\overline B(x_0;r) \subset \Omega \subset B(x_0;R),
\]
and $\Gamma$ be the family of paths in $\Omega \setminus \ov B(x_0;r)$ joining $\partial B(x_0;r)$ to $\partial \Omega$. Then
\[
M(\Ga) \geq \omega_{n-1} \Big(\log \frac{R}{r}\Big)^{1-n},
\]
where $\omega_{n-1}$ is the surface area of the unit sphere in $\mf{R}^n$.
\end{lem}
\begin{proof}
Let $\Ga_{r,R}$ be the family of paths in $B(x_0;R)\setminus \ov B(x_0;r))$ joining $\partial B(x_0;r)$ to $\partial B(x_0;R)$. Then $\Gamma < \Gamma_{r,R}$ and therefore by \eqref{mod-dec}, $M(\Ga) \geq M(\Ga_{r,R})$. On the other hand, from $\S 7.5$ of \cite{V},
\[
M(\Ga_{r,R})=\frac{\omega_{n-1}}{\big(\log (R/r)\big)^{n-1}},
\]
which completes the proof.
\end{proof}

\noindent We will only need the fact that $M(\Ga)>0$ where $\Ga$ is as above. See Lemma 7.38 of \cite{Vu} for a possibly better lower bound. Next, the following version of Hurwitz's theorem for quasiconformal mappings will be useful.

\begin{thm}\label{Hurwitz}
Let $\Omega, \Omega' \subset \mathbb{R}^n$ be bounded domains, and $x_0 \in \Omega$, $x_0' \in \Omega'$. Suppose that $f_j \in QC_K(\Om, \Om')$ is a sequence such that $f_j(x_0)=x_0'$ for all $j \ge 1$ and $\{f_j\}$ converges uniformly on compact subsets of $\Omega$ to $f: \Omega \to\mathbb{R}^n$. 

\medskip

If there exists a uniform neighbourhood $V'$ of $ x_0'$ in $\Omega'$ such that $V' \subset f_j (\Omega)$ for all large $j$, then $f$ is $K$-quasiconformal. Moreover, if $g_j=f_j^{-1}\vert_{V'}$, then, after passing to a subsequence, $\{g_j\}$ converges uniformly on compact subsets of $V'$ to a $K$-quasiconformal map $g: V' \to \Omega$ such that $f \circ g(x')=x'$ for all $x' \in V'$. In particular, $V' \subset f(\Omega)$.
\end{thm}
\begin{proof}
We know that $f$ is constant or $K$-quasiconformal. If possible, assume that $f$ is constant. Then $f \equiv x_0'$ as $f(x_0)=x_0'$. Let $\ov B(x_0;\de) \subset \Om$ and
\[
m_j'= \sup \big\{ |f_j(x) - x'_0|: x \in \ov B(x_0;\de) \big\}.
\]
Since $f_j \to x_0'$ uniformly on $\ov B(x_0;\de)$, $m_j' \to 0$. Fix $B(x_0';r') \subset V'$. Then, 
\begin{equation}\label{containment}
f_j\big(\ov B(x_0;\de)\big) \subset \ov B(x_0';m_j') \subset B(x_0';r') \subset f_j(\Om),
\end{equation}
for all large $j$. Without loss of generality, we assume that this holds for each $j \geq 1$.

\medskip

Let $\Ga$ be the family of all paths in $\Omega \setminus \ov B(x_0;\de)$ joining $\pa B(x_0;\de)$ to $\pa \Om$. Since $\Omega$ is bounded, 
\[
M(\Ga)>0
\]
by Lemma\ref{positive-capacity}. For each $j \geq 1$, let $\Ga_j=\{f_j \circ \ga: \ga \in \Ga\}$ be the image path family in $\Omega'$. Since $f_j$ is $K$-quasiconformal,
\begin{equation}\label{M-M_j}
\frac{M(\Ga_j)}{M(\Ga)} \geq \frac{1}{K}.
\end{equation}

For each $j \geq 1$, let $\Ga_j'$ be the family of paths in $B(x_0';r')\setminus \ov B(x_0';m_j')$ joining $\pa B(x_0';m_j')$ to $\pa B(x_0';r')$. In view of \eqref{containment}, each path in $\Ga_j$ has a subpath which belongs to $\Ga_j'$ and so $\Ga_j$ is minorized by $\Ga_j'$. By Theorem~6.4 of \cite{V},
\begin{equation}\label{M_j-ring}
M(\Ga_j) \leq M(\Ga_j')=\om_{n-1}\left(\log\frac{r'}{m_j'}\right)^{1-n},
\end{equation}
where the last equality follows from 7.5 of \cite{V}.

\medskip

Thus 
\[
\frac{M(\Ga)}{K} \leq \om_{n-1}\left(\log\frac{r'}{m_j'}\right)^{1-n},
\]
from \eqref{M-M_j} and \eqref{M_j-ring}. Now, observe that the right-hand side above converges to $0$, which gives $M(\Ga)=0$, a contradiction. Thus, $f$ is $K$-quasiconformal.

\medskip

Note that $g_j: V' \to \Omega$ is $K$-quasiconformal and by passing to a subsequence, $\{g_j\}$ converges uniformly on compact subsets of $V'$ to a map $g: V' \to \overline \Omega$, which is constant or $K$-quasiconformal. Since $g$ is continuous at $x_0'$ and $g(x_0')=x_0 \in \Omega$, there exists a neighbourhood $W'\subset V'$ of $x_0'$ such that $g(W') \subset \Omega$. Note that for $x' \in W'$,
\[
f_j\big(g_j(x')\big)=x'.
\]
Since $g_j(x') \to g(x') \in g(W') \subset \Omega$ and $f_j \to f$ uniformly on compact subsets of $\Omega$, we must have $f_j(g_j(x')) \to f(g(x'))$. It follows that
\[
f\big(g(x')\big)=x'.
\]
This shows that $g$ is non-constant, and hence $g$ is $K$-quasiconformal on $V$. In particular, $g: V' \to \Omega$. The same argument as above now shows that $f(g(x'))=x'$ for all $x' \in V'$. This in turn implies that
\[
V'=f\big(g(V')\big) \subset f(\Omega),
\]
which completes the proof.
\end{proof}

The same arguments also give us the following:

\begin{thm}\label{Hurwitz-2}
Let $\{\Om_j\}$ and $\{\Om_j'\}$ be sequences of uniformly bounded domains in $\mathbb{R}^n$ that converge in the local Hausdorff sense to bounded domains $\Om$ and $\Om'$, respectively. Let $x_0\in \Omega$ and $x_0' \in \Omega'$. Suppose $f_j \in QC_K(\Om_j, \Om'_j)$ is a sequence such that $f_j(x_0)=x_0'$ for all large $j$, and $\{f_j\}$ converges uniformly on compact subsets of $\Omega$ to $f: \Omega \to \mathbb{R}^n$. 

\medskip

If there exists a uniform neighbourhood $V'$ of $ x_0'$ in $\Omega'$ such that $V' \subset f_j (\Omega)$ for all large $j$, then $f$ is $K$-quasiconformal. Moreover, if $g_j=f_j^{-1}\vert_{V'}$, then, after passing to a subsequence, $\{g_j\}$ converges uniformly on compact subsets of $V'$ to a $K$-quasiconformal map $g: V' \to \Omega$ such that $f \circ g(x')=x'$ for all $x' \in V'$. In particular, $V' \subset f(\Omega)$. 
\end{thm}

\begin{thm}\label{extremal}
Let $\Om \subset \mathbb{R}^n$ be a bounded domain and $x_0 \in \Omega$. Then there exists a $K$-quasiconformal map $f: \Omega \to \mathbb{B}^n$ such that $f(x_0)=0$ and $B(0;s_{\Omega}^K(x_0)) \subset f(\Omega)$.
\end{thm}
\begin{proof}
Choose an increasing sequence $r_j \ra s_{\Om}^K(x_0)$ such that for each $j$ there exists an $f_j \in QC_K(\Om, \mf{B}^n)$ such that $f_j(x_0)=0$ and $B(0;r_j) \subset f_j(\Omega)$.

\medskip

The family $\{f_j:j \ge 1\}$ is normal and therefore, after passing to a subsequence, $\{f_j\}$ converges uniformly on compact subsets of $\Om$ to $f: \Om \to \overline{\mathbb{B}^n}$ which is constant or $K$-quasiconformal. Since $\{r_j\}$ is increasing, $B(0;r_1) \subset f_j(\Om)$ for all $j \geq 1$. Therefore, by Theorem~\ref{Hurwitz}, $f$ is $K$-quasiconformal.

\medskip

To prove that $B(0;s_{\Om}^K(z_0)) \subset f(\Om)$, it suffices to prove that $B(0;r_j) \subset f(\Om)$ for each $j$. Fix $j$. Then for each $i>j$, $B(0;r_j) \subset B(0;r_i) \subset f_i(\Om)$. By Theorem~\ref{Hurwitz} again, $B(0;r_j) \subset f(\Om)$ as required.
\end{proof}

\begin{prop}\label{s-ball}
Let $\Omega \subset \mathbb{R}^n$ be a bounded domain and $K\ge 1$. If there exists $x_0 \in \Omega$ such that $s^K_{\Omega}(x_0)=1$, then $\Om$ is $K$-quasiconformally equivalent to $\mathbb{B}^n$, and hence $s^K_{\Omega}(x)=1$ for all $x \in \mathbb{B}^n$.
\end{prop}
\begin{proof}
Let $f$ be an extremal map at $x_0$. Then
\[
\mf{B}^n \subset f(\Om) \subset \mf{B}^n
\]
which implies that $f(\Om)=\mf{B}^n$ and hence $f$ is a $K$-quasiconformal homeomorphism of $\Om$ onto $\mathbb{B}^n$.

\medskip

For the second part, given $x \in \Om$, let $T$ be a conformal self-map of $\mathbb{B}^n$ with $T(f(x))=0$. Then $T\circ f: \Omega \to \mathbb{B}^n$ is a $K$-quasiconformal homeomorphism such that $T(f(x))=0$ and
\[
\mathbb{B}^n  \subset T\circ f(\mathbb{B}^n).
\]
Therefore $s^K_{\Omega}(x) \geq 1$ and hence  $s^K_{\Omega}(x) = 1$.

\end{proof}

\begin{prop}
Let $\Om \subset \mf{R}^n$ be a bounded domain and $x \in \Om$. Then
\[
\sup_{K \in [1,\infty)}s^K_{\Om}(x)=S_{\Om}(x).
\]
\end{prop}
\begin{proof}
Let $f$ be an extremal for $s^K_{\Omega}(x)$, i.e., $f: \Om \to \mf{B}^n$ is a $K$-quasiconformal 
map such that $f(x)=0$ and $B(0;s^K_{\Om}(x)) \subset f(\Om)$. Therefore $s^K_{\Om}(x)$ is a candidate for $S_{\Om}(x)$ and hence $s^K_{\Om}(x) \leq S_{\Om}(x)$. Taking the supremum over $K$,
\[
\sup_{K \in [1,\infty)} s_{\Om}^K(x) \leq S_{\Om}(x).
\]
To obtain the reverse inequality, observe that if $r$ is a candidate for $S_{\Om}(x)$, then there exists $K\in [1,\infty)$ such that $r$ is also a candidate for $s_{\Om}^K(x)$. Therefore
\[
r \leq s_{\Om}^K(x) \leq \sup_{K\in [1,\infty)}s_{\Om}^K(x).
\]
Taking the supremum over $r$, we obtain
\[
S_{\Om}(x) \leq \sup_{K\in [1,\infty)}s_{\Om}^K(x)
\]
and this completes the proof.
\end{proof}

\begin{prop}\label{s-dec}
Let $\Om \subset \mf{R}^n$ be a bounded domain that is not quasiconformally equivalent to $\mf{B}^n$. Fix $x \in \Om$. If $1 \leq K_1 < K_2 $, then $s^{K_1}_{\Om}(x)< s^{K_2}_{\Om}(x)$.
\end{prop}
\begin{proof}
Let $ r =s^{K_1}_{\Om}(x)$ and let $f: \Om \to \mf{B}$ be a $K_1$-quasiconformal map such that $f(x)=0$ and $B(0;r) \subset f(\Om)$. Let $0<a<1$ and $h: \mf{B}^n \to \mf{B}^n$ be the radial map defined by
\[
h(x)=
\begin{cases}
\vert x\vert^{a-1}x & \text{if $x \neq 0$},\\
0 & \text{if $x=0$}.
\end{cases}
\]
Then $h$ is a $\vert a\vert^{1-n}$-quasiconformal homeomorphism of $\mf{B}^n$ onto itself (see \cite[16.2]{V}). Taking
\[
a=\left(\frac{K_2}{K_1}\right)^{\frac{1}{1-n}}<1,
\]
$h$ is a $K_2/K_1$-quasiconformal homeomorphism of $\mf{B}^n$ onto itself. Therefore, $h \circ f: \Om \to \mf{B}^n$ is a $K_2$-quasiconformal map and
\[
B(0;r^{a})=h\big(B(0;r)\big) \subset h\circ f(\Om),
\]
which shows that $s^{K_2}_{\Om}(x)\geq r^a$. Since $0<r<1$ and $0<a<1$, we have $r^a>r$. It follows that
\[
s^{K_2}_{\Om}(x) \geq r^{a}>r=s^{K_1}_{\Om}(x).
\]
This completes the proof.
\end{proof}

\subsection{Continuity of $s^K_{\Om}$}


The proof of the following continuity result is analogous to that for the usual squeezing function \cite[Theorem~3.1]{DGZ2012} except that instead of the invariant Kobayashi metric we use the quasi-invariant $\mu$-metric.

\begin{thm}
Let $\Omega \subset \mf{R}^n$ be a bounded domain. Then $s^K_{\Omega}$ is continuous on $\Omega$.
\end{thm}
\begin{proof}
Let $x_0\in \Omega$ and choose $x_j \ra x_0$. The proof has two parts.

\medskip

\noindent {\sf Step I}: To see that $\liminf_{j \to \infty} s^K_{\Omega}(x_j) \geq s^K_{\Omega}(x_0)$, let $\epsilon>0$ and let $f: \Omega \to \mathbb{B}^n$ be a $K$-quasiconformal map such that $f(x_0)=0$ and $B(0;s^K_{\Omega}(x_0))\subset f(\Omega)$. Since $f$ is continuous, there exists an integer $N$ such that $\vert f(x_j)\vert <\ep$ for $j \geq N$. Define $f_j: \Omega \to \mf{R}^n$ by
\[
f_j(x)=\frac{f(x)-f(x_j)}{1+\ep},
\]
for $j\geq N$. Then $f_j$ is $K$-quasiconformal (since it is the composition of a $K$-quasiconformal map with a conformal map), $f_j: \Om \to \mathbb{B}^n$, $f_j(x_j)=0$, and
\[
B\left(0;\frac{s^K_{\Omega}(x_0)-\ep}{1+\ep}\right) \subset f_j(\Omega).
\]
To see the above inclusion, note that, if 
\[
\vert y \vert < \frac{s^K_{\Omega}(x_0)-\ep}{1+\ep}
\]
then
\[
 \big\vert (1+\ep)y+f(x_j)\big\vert \leq  \big\vert (1+\ep)y\big\vert +\ep < s^K_{\Omega}(x_0),
\]
and therefore
\[
(1+\ep)y+f(x_j) \subset B\big(0;s^K_{\Omega}(x_0)\big) \subset f(\Omega),
\]
which shows that $y \in f_j(\Omega)$, as required. Thus,
\[
s^K_{\Omega}(x_j) \geq \frac{s^K_{\Omega}(x_0)-\ep}{1+\ep}
\]
and hence
\[
\liminf_{j \to \infty} s^K_{\Omega}(x_j) \geq s^K_{\Omega}(x_0).
\]

\medskip

\noindent {\sf Step II}: To see that $\limsup_{j \to \infty}s^K_{\Omega}(x_j) \leq s^K_{\Omega}(x_0)$, first note that $\mu_{\Omega}(x_j,x_0) \to 0$.

\medskip

For each $j$, let $f_j: \Omega \to \mathbb{B}^n$ be a $K$-quasiconformal mapping such that $f_j(x_j)=0$ and $B(0;s^K_{\Omega}(x_j)) \subset f_j(\Omega)$. Now
\[
\mu_{\mathbb{B}^n}\big(0,f_j(x_0)\big)=\mu_{\mathbb{B}^n}\big(f_j(x_j),f_j(x_0)\big) \leq  \mu_{f_j(\Omega)}\big(f_j(x_j),f_j(x_0)\big) \leq K \mu_{\Omega}(x_j,x_0) \to 0.
\]
Since the topology induced by $\mu_{\mf{B}^n}$ is the same as the Euclidean topology, we conclude that $f_j(x_0) \to 0$. Let $\epsilon>0$. Then, there exist $N \ge 1$ such that $\vert f_j(x_0)\vert <\epsilon$ for $j \geq N$.

\medskip

Note that $F_j: \Omega \to \mathbb{R}^n$ defined as
\[
F_j(x)=\frac{f_j(x)-f_j(x_0)}{1+\epsilon}
\]
satisfies $F_j: \Omega \to \mathbb{B}^n$, $F_j(x_0)=0$, and
\[
B\left(0;\frac{s^K_{\Omega}(x_j)-\epsilon}{1+\epsilon}\right) \subset F_j(\Omega).
\]
To see the above inclusion, note that, if
\[
\vert y \vert < \frac{s^K_{\Omega}(x_j)-\epsilon}{1+\epsilon} 
\]
then
\[
\big\vert (1+\epsilon)y + f_j(x_0) \big\vert < \big\vert (1+\epsilon)y \big\vert +\epsilon < s^K_{\Omega}(x_j),
\]
and therefore
\[
(1+\epsilon)y + f_j(x_0) \subset B\big(0;s^K_{\Omega}(x_j)\big) \subset f_j(\Omega) \, , 
\]
which shows that $y \in F_j(\Omega)$, as required. Thus
\[
s^K_{\Omega}(x_0) \geq \frac{s^K_{\Omega}(x_j)-\epsilon}{1+\epsilon}
\]
and hence
\[
s^K_{\Omega}(x_0) \geq \limsup_{j \to \infty} s^K_{\Omega}(x_j).
\]
This establishes the continuity of $s^K_{\Omega}$.
\end{proof}

\subsection{Dependence of $s^K_{\Om}$ on $\Om$} 

We investigate the relationship between the $K$-squeezing functions of a monotone sequence of domains and that of its limit domain. We begin with the case of an increasing sequence of domains.

\begin{thm}
Let $\Omega \subset \mathbb{R}^n$ be a bounded domain and $\Omega_j \subset \Om_{j+1}$ be an increasing sequence of sub-domains of $\Omega$ such that $\cup_{j} \Omega_j=\Omega$. Let $K \ge 1$. Then
\[
s^K_{\Omega_j}(x) \to s^K_{\Omega}(x).
\]
for every $x \in \Omega$.
\end{thm}
\begin{proof}
Fix $x_0 \in \Omega$. The proof is divided into two steps.

\medskip

\noindent {\sf Step I}: To see that $s^K_{\Omega}(x_0) \geq \limsup_{j \to \infty}s_{\Omega_j}(x_0)$, let $f_j: \Omega_j \to \mathbb{B}^n$ be a $K$-quasiconformal map such that $f_j(x_0)=0$ and $B(0;s^K_{\Omega_j}(x_0)) \subset f_j(\Omega_j)$. Passing to a subsequence, $\{f_j\}$ converges uniformly on compact subsets of $\Omega$ to $f: \Omega \to \overline{\mathbb{B}^n}$.

\medskip

Fix $j_0 \ge 1$ and note that for $j \geq j_0$,
\begin{equation}\label{s-lowbd}
s^K_{\Omega_j}(x_0)\geq \frac{d(x_0,\pa \Omega_j)}{\text{diam}(\Omega_j)} \geq \frac{d(x_0,\Omega_{j_0})}{\text{diam}(\Omega)}>0.
\end{equation}
This shows that there exists a uniform $\de>0$ such that $B(0;\de) \subset f_j(\Omega_j)$ for all $j \geq j_0$. Hence, by Theorem~\ref{Hurwitz-2}, $f$ is $K$-quasiconformal. By invariance of domains, $f(\Omega)$ is an open subset of $\overline{\mathbb{R}^n}$ and so $f: \Omega \to \mathbb{B}^n$. Also, $f(x_0)=0$.

\medskip

Let $r=\limsup_{j \to \infty}s^K_{\Omega_j}(x_0)$. Note that $r > 0$ by \eqref{s-lowbd}. Then there is a subsequence of $(s^K_{\Omega_j}(x_0))$, say $(s^K_{\Omega_{j_k}}(x_0))$, which converges to $r$. Let $\ep\in (0,r)$. Then there exists $k_0\ge 1$ such that
\[
r-\ep< s^K_{\Omega_{j_k}}(x_0), \quad k\geq k_0
\]
and therefore
\[
B(0;r-\ep) \subset f_{j_k}(\Omega_{j_k}).
\]
By Theorem~\ref{Hurwitz-2}, $B(0;r-\ep) \subset f(\Omega)$. Since $f: \Omega \to \mathbb{B}^n$ is $K$-quasiconformal with $f(x_0)=0$, we must have $s^K_{\Omega}(x_0) \geq r-\ep$. Since $\ep \in (0,r)$ is arbitrary, this implies that $s^K_{\Omega}(x_0)\geq r$.

\medskip

\noindent {\sf Step II}: To see that $s^K_{\Omega}(x_0) \leq \liminf_{j \to \infty} s^K_{\Omega_j}(x_0)$, let $\phi: \Omega \to \mathbb{B}^n$ be a $K$-quasiconformal map such that $\phi(x_0)=0$ and $B(0;s^K_{\Omega}(x_0)) \subset \phi(\Omega)$. Let $\ep>0$. Then
\[
\ov B\big(0;s^K_{\Omega}(x_0)-\ep\big) \subset B\big(0;s^K_{\Omega}(x_0)\big)\subset \phi(\Omega).
\]
Therefore $\phi^{-1}(\ov B(0;s^K_{\Omega}(x_0)-\ep))$ is a compact subset of $\Omega$ and therefore is contained in $\Omega_j$ for all large $j$. This implies that $B(0;s^K_{\Omega}(x_0)-\ep)) \subset \phi(\Om_j)$ for all large $j$. Since $\Om_j \subset \Om$, $\phi$ is quasiconformal on $\Omega_j$. Therefore
\[
s^K_{\Omega_j}(x_0) \geq s^K_{\Omega}(x_0)-\ep
\]
for all large $j$ and hence
\[
\liminf_{j \to \infty}s_{\Omega_j}(x_0)\geq s^K_{\Omega}(x_0)-\ep.
\]
Since $\ep>0$ is arbitrary, we must have $\liminf_{j \to \infty}s^K_{\Omega_j}(x_0)\geq s^K_{\Omega}(x_0)$.
\end{proof}

\begin{thm}
Let $\Omega \subset \mathbb{R}^n$ be a bounded domain and $\Omega_j \supset \Om_{j+1}$ be a decreasing sequence domains such that $\cap_{j} \Omega_j=\Omega$. Let $K \ge 1$. Then
\[
s^K_{\Omega}(x) \geq \limsup  s^K_{\Omega_j}(x) .
\]
for every $x \in \Omega$.
\end{thm}
\begin{proof}
The proof is essentially Step I of the previous theorem. We only need to modify \eqref{s-lowbd} as follows:
\begin{equation}
s^K_{\Omega_j}(x_0)\geq \frac{d(x_0,\pa \Omega_j)}{\text{diam}(\Omega_j)} \geq \frac{d(x_0,\Omega)}{\text{diam}(\Omega_1)}>0.
\end{equation}
for all $j \geq 1$.
\end{proof}

\noindent As a concluding remark in this section, the $K$-squeezing function $s^K_{\Om}$ allows us to interpret a result of Gehring--V\"ais\"al\"a \cite{GV} on the class of domains in $\mbb R^3$ that are quasiconformally equivalent to the ball $\mbb B^3$. To explain this, let $\Om \subset \mbb R^3$ be homeomorphic to the ball $\mbb B^3$. Following \cite{GV}, let
\[
K(\Om)=\inf\big\{K \in [1,\infty): \text{$\Om$ is $K$-quasiconformally equivalent to $\mathbb{B}^3$}\big\}
\]
and let $\Xi$ denote the class of domains $\Om \subset \mbb R^3$ with $K(\Om) < \infty$. Let $\mathcal B$ be the class of bounded domains in $\mbb R^3$. The distance between $\Om, \Om' \subset \Xi$ is defined as
\[
d(\Om, \Om') = \inf \log K 
\]
where the infimum is taken over all $K \in [1, \infty)$ such that there exists a $K$-quasiconformal homeomorphism from $ \Om $ onto $ \Om' $. A pair of domains $\Om, \Om'$ with $d(\Om, \Om') = 0$ were also identified in \cite{GV}. In this setup, it is shown in $\S 11$ of \cite{GV} that $(\Xi, d)$ is a complete, nonseparable metric space.

\medskip

Proposition \ref{s-ball} shows that for a bounded domain $\Om$, $s^K_{\Om} \equiv 1$ is equivalent to $\Om$ being $K$-quasiconformally equivalent to the ball $\mbb B^n$. Thus, by viewing $s^K$ as a domain functional defined on $\mathcal B$ with values in the interval $(0, 1]$, this means that the fibre $(s^K_{\Om})^{-1}(1)$ consists precisely of those domains in $\mathcal B$ that are $K$-quasiconformally equivalent to $\mbb B^n$. The union of these fibres, as $K$ varies over $[1, \infty)$, is exactly the intersection $\Xi \cap \mathcal B$. Thus, the result of Gehring--V\"ais\"al\"a can be understood as saying that 
\[
\Xi \cap\mathcal B = \bigcup_{K \ge 1} (s^K_{\Om})^{-1}(1)
\]
with the distance $d$, defined above, is complete and nonseparable when $n=3$.

\subsection{Boundary behaviour of $s^K_{\Om}$} We will now look at the asymptotic behaviour of $s^K_{\Om}$ near $\pa \Om$. This is done first for planar domains and then when the domain $\Om \subset \mbb R^n$, $n \ge 3$. As will be seen, the planar case is a direct consequence of what is known for the usual (holomorphic) squeezing function. When $\Om \subset \mbb R^n$, $n \ge 3$, we use a conformal scaling technique to understand the boundary behaviour of $s^K_{\Om}$ near {\it spherically extreme} points on $\pa \Om$.

\begin{thm}
Let $\Omega \subset \mathbb{C}$ be a finitely connected domain such that no boundary component of $\Omega$ is a singleton, and fix $K \ge 1$. Then
\[
\lim_{z \to \pa \Omega}s^K_{\Omega}(z)=1.
\]
\end{thm}
\begin{proof}
In complex dimension one, conformal embeddings are precisely holomorphic embeddings, and so $s^1_{\Omega}$ coincides with the usual squeezing function of $\Omega$. Therefore, by \cite[Theorem~5.2]{DGZ2012},
\[
\lim_{z \to \pa \Omega} s^1_{\Omega}(z)=1.
\]
Since $s^1_{\Omega} \leq s^K_{\Omega}\leq 1$, the proof of the theorem follows.
\end{proof}

Let $\Omega \subset \mathbb{R}^n$ be a bounded domain, $p \in \partial\Omega$, and suppose that $\partial \Omega$ is smooth near $p$. Write $x = (x_1, x') \in \mbb R \times \mbb R^{n-1}$, where $x' = (x_2, \ldots, x_n)$. Following \cite[Section $3$]{KZ2016}, we say that $p$ is a {\it spherically extreme} boundary point of $\Omega$ if there exists a ball $B \subset \mathbb{R}^n$ such that
\begin{equation}\label{sp-ext}
\Omega \subset B \quad \text{and} \quad p \in \partial \Omega \cap \partial B.
\end{equation}

We first show that the property of being spherically extreme is open, i.e., if $p \in \partial \Omega$ is spherically extreme, then every boundary point of $\Omega$ sufficiently close to $p$ is also spherically extreme. In the holomorphic case, this was done in \cite{KZ2016}. Since the notion of a spherically extreme boundary point is invariant under rigid motions of $\mf{R}^n$, we can assume without loss of generality that $p=0$ and the outward unit normal to $\partial \Omega$ at the origin is along the $(-1,0')$ direction. Suppose that near $x=0$, the graph of $\partial \Omega$ is given by $x_1=\phi(x')$. Then $\phi(0')=0$ and $\nabla\phi(0')=0$. Thus, the Taylor expansion of $\phi$ at $x=0$ is given by
\begin{equation}\label{phi-TE}
\phi(x')=\sum_{i,j=2}^n \phi_{ij}(0')x_ix_j+O(\vert x'\vert^3).
\end{equation}
Suppose $B$ is a ball of radius $R$ satisfying \eqref{sp-ext}. Since both $\Om$ and $\mbb B^n$ are smooth, $\partial B$ and $\partial \Omega$ must have the same outward unit normal at $0$, and, therefore $B=B((R,0');R)$. Thus, the defining equation of $\partial B$ is
\[
(x_1-R)^2 +\vert x'\vert^2=R^2.
\]
Near $x=0$, the above sphere is the graph of $x_1=\psi(x')$ where
\begin{equation}\label{psi-TE}
\psi(x')=R-\sqrt{R^2-\vert x'\vert^2} = R\left(1-\sqrt{1-\frac{\vert x'\vert^2}{R^2}}\right) = \frac{\vert x'\vert^2}{2R}+ \frac{\vert x'\vert^4}{8R^3}+O(\vert x'\vert^6).
\end{equation}
Since $\Omega \subset B$, $\partial \Omega$ lies to the right of $\pa B$ near $x=0$, and hence $\phi(x')\geq \psi(x')$ for $x'$ near $0$. From \eqref{phi-TE} and \eqref{psi-TE}, it follows that
\begin{equation}\label{sp-ext-ch}
H_\phi(0') \geq \frac{1}{2R}I,
\end{equation}
where $H_{\phi}(x')$ denotes the Hessian matrix of the function $x_1=\phi(x')$ at $x'$ and $I$ is the $(n-1) \times (n-1)$ identity matrix. Here the matrix inequality $A \geq B$ (resp. $A>B$) means that the matrix $A-B$ is semi-positive definite (resp. positive definite).

\begin{lem}\label{sp-ext-ch-lem}
Let $S>R$. Then for every $\zeta \in \partial \Omega$ sufficiently close to $p$, there exists a ball $B_{\zeta}$ of radius $S$ such that $\Omega \subset B_\zeta$ and $\zeta \in \partial \Omega \cap \partial B_\zeta$.
\end{lem}
\begin{proof}
Let $\zeta\in \partial \Omega$. Denote the outward unit normal to $\partial \Omega$ at $\zeta$ by $n_\zeta$.  Let $c_\zeta=\zeta-Sn_\zeta$ and $B_\zeta=B(c_\zeta;S)$. Evidently $\zeta \in \partial \Omega \cap \partial B_\zeta$. We will show that $\Omega \subset B_\zeta$ if $\zeta$ is sufficiently close to $p$ in the following two steps:

\medskip

\noindent {\sf Step I}: There exists $\ep>0$ such that $\Omega \cap B(p;\ep) \subset B_\zeta$ for all $\zeta\in \partial \Omega$ sufficiently close to $p$. To see this, choose $\alpha>0$ such that $S>\alpha>R$. Then from \eqref{sp-ext-ch} and by continuity of $H_{\phi}$,
\begin{equation}\label{H-0}
H_{\phi}(x') > \frac{1}{2\alpha}I,
\end{equation}
for all $x'$ near $0'$.

\medskip

Note that $\partial B_\zeta$ near $\zeta$ is given by the graph of the function
\[
x_1=\psi_\zeta(x')=(c_\zeta)_1-\sqrt{S^2-\vert x'-c_\zeta'\vert^2}.
\]
Since $\partial \Omega$ is smooth near $p$, $c_{\zeta}$ and hence the function $\psi_\zeta$, varies smoothly in $\zeta$ near $p$. Also, at $p$, by expanding $\psi_p$ as in \eqref{psi-TE}, we have
\[
H_{\psi_{p}}(0')=\frac{1}{2S}I.
\]
As $1/S < 1/\alpha$, by continuity of $H_{\psi_{\zeta}}(x')$ in $\zeta$ and $x'$,
\begin{equation}\label{H-q}
H_{\psi_{\zeta}}(x')< \frac{1}{2\alpha}I,
\end{equation}
for $\zeta$ near $p$ and $x'$ near $0'$. Combining \eqref{H-0} and \eqref{H-q}, we see that
\[
H_{\phi}(x') > H_{\psi_{\zeta}}(x'),
\]
for $x'$ near $0'$. This implies that if $\zeta$ is sufficiently close $p$, then the graph of $x_1=\phi(x')$ lies to the right of the graph of $x_1=\psi_\zeta(x')$ near $p$. In other words, there exists $\epsilon>0$ such that
\[
\Omega \cap B(p;\epsilon) \subset B_\zeta,
\]
if $\zeta \in \partial \Omega$ is sufficiently close to $p$.

\medskip

\noindent {\sf Step II}: Now we will show that $\overline \Omega \setminus B(p;\epsilon) \subset B_\zeta$, i.e.,
\[
\vert x-c_\zeta\vert^2 <S^2, \quad x \in \overline \Omega \setminus B(p;\epsilon),
\]
for every $\zeta\in \partial \Omega$ sufficiently close to $p$. Since $\overline\Omega \setminus B(p;\epsilon)$ is compact and the function $\vert x -c_\zeta\vert^2$ is continuous in $(x,\zeta)$ for $x \in \overline \Omega \setminus B(p;\epsilon)$ and $\zeta \in \partial \Omega$ near $p$, it is enough to show that
\begin{equation}\label{ineq-al}
\vert x-c_p\vert^2< \beta, \quad x \in \overline \Omega \setminus B(p;\epsilon),
\end{equation}
for some $\beta<S^2$.

\medskip
 
To prove this, recall that $B=B((R,0');R)$ and $B_p=B((S,0');S)$. As $R<S$, we have $B\subset B_p$ with $p=0$ being the only common boundary point. Thus,
\[
\overline \Omega  \subset \overline B  \subset B_p \cup \{p\},
\]
and so
\[
\vert x-c_p\vert^2 < S^2, \quad x \in \overline \Omega \setminus B(p;\ep).
\]
By continuity of the function $\vert x -c_p\vert^2$ on the compact set $\overline \Omega \setminus B(p;\ep)$, it follows that there exists $\be<S^2$ such that \eqref{ineq-al} holds. This completes Step II and the proof of the lemma.
\end{proof}

The boundary behaviour of $s^K_{\Om}$ near spherically extreme boundary points is addressed below. At first, we look at the special case when such boundary points are approached along the inner normal. Next, the general case reduces to this special case by using a set of suitable conformal transformations.

\begin{thm}\label{bd-bev-s}
Let $\Omega \subset \mathbb{R}^n$ be a bounded domain and 
$\{p_j\}$ be a sequence in $\Omega$ that converges to $p \in \partial \Omega$. Assume that $ \partial \Omega$ is smooth near $ p$ and $p$ is spherically extreme. Then $s^K_{\Omega}(p_j) \to 1$.
\end{thm}
\begin{proof}
By means of a translation and an orthogonal transformation, which preserve the squeezing function $s^K_{\Om}$, we can assume that $p=0$ and that the outward normal to $\partial \Omega$ at the origin points in the direction $(-1,0')$. There are two cases to consider.

\medskip

\noindent {\sf Case I}: The sequence $\{p_j\}$ approaches $p$ along the inward normal.

\medskip

Let $\delta_j=d(p_j, \partial \Omega)$. Then $p_j=(\delta_j,0')$. Since $p=0$ is a spherically extreme boundary point of $\Omega$, there exists $R>0$ such that
\[
\Omega \subset B=B\big((R,0'); R\big).
\]
The conformal affine maps $T_j: \mathbb{R}^n \to \mathbb{R}^n$ defined by
\[
T_j(x)=\frac{1}{\de_j}(x-p_j)
\]
satisfy $T_j(p_j) = 0$ for all $j \ge 1$. Furthermore, 
\[
T_j(B)=B\left(\left(\frac{R}{\delta_j}-1,0'\right); \frac{R}{\delta_j}\right) \subset \{x_1>-1\},
\]
since $(-1,0) \in T_j(B)$, and the tangent-plane to $T_j(B)$ at $(-1,0)$ is $\{x: x_1=-1\}$. Therefore
\[
T_j(\Omega) \subset T_j(B) \subset \{x_1>-1\},
\]
for each $j \geq 1$. Write $\Om_j = T_j(\Omega)$ and $\mathcal{H} = \{x: x_1>-1\}$, which is a half-space. By conformal invariance
\[
s^1_{\Om}(p_j)= s^1_{\Om_j}(0),
\]
and hence, it is enough to establish that $s^1_{\Omega_j}(0)  \to 1$. In fact, as $s^1_{\Om_j}(0) \leq 1$, it is enough to show that $1\leq \liminf_{j\to \infty}s^1_{\Omega_j}(0)$.

\medskip

The domains $\Omega_j \to \mathcal{H}$ in the local Hausdorff sense. In particular, if $K \subset \mathcal{H}$ is a compact set, then $K \subset \Omega_j$ for all large $j$.

\medskip

Let $\phi: \mathcal{H} \to \mathbb{B}^n$ be a conformal homeomorphism such that $\phi(0)=0$. Let $\ep>0$. Then $\phi^{-1}(\ov B(0;1-\ep))$ is compactly contained in $\mathcal{H}$ and is therefore contained in $\Omega_j$ for all large $j$. Since $\Omega_j \subset \mathcal{H}$, $\phi$ is a $1$-quasiconformal map of $\Om_j$ into $\mathbb{B}^n$. Thus, $B(0;1-\ep) \subset \phi(\Om_j)$ and therefore, the map $\phi$ is a candidate for $s^1_{\Omega_j}(0)$ 
for all large $j$. Therefore
\[
1-\ep \leq s^1_{\Omega_j}(0),
\]
for all large $j$. It follows that
\[
1-\ep \leq \liminf_{j \to \infty}s^1_{\Omega_j}(0),
\]
and hence
\[
1\leq \liminf_{j \to \infty}s^1_{\Omega_j}(0),
\]
as required. Since $s^1_{\Om} \le s^K_{\Om}$, for $K \ge 1$, it follows that $s^K_{\Om}(p_j) \ra 1$.

\medskip

\noindent {\sf Case II}: The sequence $\{p_j\}$ approaches $p$ in an unrestricted manner.

\medskip

In this case, choose positive constants $r,R$ such that for each $\zeta \in \partial \Omega \cap B(0;r)$, there exists a ball $B_\zeta$ of radius $R$ such that
\[
\Omega \subset B_{\zeta} \quad \text{and} \quad \zeta \in \partial \Omega \cap \partial B_{\zeta}.
\]
Without loss of generality, we assume that each $p_j \in B(0;r)$ and each $p_j$ admits a point $q_j \in \partial \Omega \cap B(p;r)$ which is closest to it. Thus, $d(p_j,\partial \Omega)=\vert q_j-p_j\vert =\de_j > 0$.

\medskip

Let $O_j: \mathbb{R}^n \to \mathbb{R}^n$ be an orientation preserving orthogonal transformation that maps the unit vector $n_j=(q_j-p_j)/\vert q_j-p_j\vert$ (which is the outward unit normal to $\partial \Omega$ at $q_j$) to $(-1,0')$. The affine conformal map $A_j: \mathbb{R}^n \to \mathbb{R}^n$
\[
A_jx=O_j(x-q_j)
\]
satisfies $A_j(q_j)=0, A_j(p_j) = (\de_j,0')$, and $A_j \to {Id}_n$.

\medskip

Recall that $\Omega \subset B_{q_j}$ and $q_j \in \partial \Omega_j \cap \partial B_{q_j}$. Therefore $B_{q_j}$ and $\pa \Omega$ have the same outward unit normal at $q_j$, which is $n_j$. It follows that $B_{q_j}=B(q_j-Rn_j;R)$.  Hence
\[
A_j(\Omega) \subset A_j(B_{q_j}) = B\big((R,0');R\big)=B,
\]
say.

\medskip

Thus we have a sequence of domains $A_j(\Omega)$ contained in a uniform ball $B$, sharing a common boundary point $p=0$, and a common inner normal at $0$, and a sequence of points $A_j(p_j)=(\delta_j,0)$ on this common inner normal that converges to $0$. The final step is to compute
\[
s^1_\Omega(p_j)=s^1_{A_j(\Omega)}\big(A_j(p_j)\big).
\]
Thus we are in the same situation as in the previous theorem, except that on the right-hand side above, the domains are varying. This can be handled exactly as in the previous case.

\medskip

Consider the affine conformal maps $T_j: \mathbb{R}^n \to \mathbb{R}^n$
\[
T_j(x)=\frac{1}{\de_j}\big(x-A_j(p_j)\big),
\]
and set $\Omega_j=T_j(A_j(\Omega))$. Then
\[
s^1_{A_j(\Omega)}\big(A_j(p_j)\big) = s^1_{\Omega_j}(0),
\]
and hence it is enough to establish that $s^1_{\Omega_j}(0)  \to 1$. In fact, as $s^1_{\Om_j}(0) \leq 1$, it is enough to prove that $1\leq \liminf_{j\to \infty}s^1_{\Omega_j}(0)$.

\medskip

Observe that
\[
\Omega_j=T_j\big(A_j(\Omega)\big) \subset T_j(B)=B\left(\left(\frac{R}{\delta_j}-1,0'\right); \frac{R}{\delta_j}\right) \subset \{x_1>-1\}=\mathcal{H} \, ,
\]
say. Moreover, $\Om_j \to \mathcal{H}$ in the Hausdorff sense. To see this, let $\psi(x)=-d(x,\partial \Omega)$ which is a smooth local defining function for $\Omega$ at $0$ and satisfies $\vert \nabla\psi(\zeta)\vert =1$ for $\zeta \in \pa \Om$ near $0$. Without loss of generality, assume that $\psi$ is defined and smooth on $B(0;r)$. Then
\[
\psi_j=\frac{1}{\de_j}(\psi \circ A_j^{-1} \circ T_j^{-1})
\]
is a local defining function for $\Om_j$ near $x=0$. Let $\psi_{\infty}(x)=-1-x_1$.Then, the above Hausdorff convergence follows if we show that
\[
\psi_j \to \psi_{\infty}
\]
uniformly on compact subsets of $\mathbb{R}^n$. We do this in the following two steps:

\medskip

\noindent {\sf Step I}: $\psi\circ A_j^{-1} \to \psi$, together with all derivatives, uniformly on compact subsets of a neighbourhood of $x=0$. To see this, observe that $A_j^{-1} \to Id_n$ and so by taking $r_1<r/2$, $A_j^{-1}$ maps $B(0;r_1)$ into $B(0;r)$ for all large $j$. Therefore $\psi \circ A_j^{-1}$ is defined on $B(0;r_1)$ and the convergence of $\psi \circ A_j^{-1}$ to $\psi$ is now evident.

\medskip

\noindent {\sf Step II}: $\psi_j \to \psi_{\infty}$ uniformly on compact subsets of $\mathbb{R}^n$. To see this, first note that
\[
T_j^{-1}(x)=A_j(p_j)+\de_j x \to 0
\]
uniformly on compact subsets of $\mathbb{R}^n$. Therefore if $K$ is a compact subset of $\mathbb{R}^n$, then for all large $j$, $T_j^{-1}(K)$ lies in $B(0;r_1)$ and therefore $\psi \circ A_j^{-1} \circ T_j^{-1} $ is defined on $K$. Observe that
\begin{align*}
(\psi \circ A_j^{-1} \circ T_j^{-1})(x) & = \psi \circ A_j^{-1} \big(A_j(p_j)+\de_j x\big)\\
&= \psi (p_j) + \left\langle \nabla (\psi \circ A_j^{-1})\big(A_j(p_j)\big), (\de_j x) \right\rangle + O(\de_j^2).
\end{align*}
Since $\psi(p_j)=-\de_j$, dividing both sides by $\de_j$,
\[
\psi_j(x)=-1+ \left\langle \nabla (\psi \circ A_j^{-1})\big(A_j(p_j)\big), x \right\rangle + O(\de_j)
\]
which converges uniformly on compact subsets of $\mf{R}^n$ to
\[
-1+\nabla \psi(0) \cdot x=-1-x_1 = \psi_{\infty}(x)
\]
as $\nabla \psi(0)=(-1,0')$ by our choice of $\psi$. It follows that $\Om_j \to \Om$ in the local Hausdorff sense.

\medskip

To complete the proof of the theorem, we now show that $s^1_{\Omega_j}(0) \to 1$. Let $\phi: \mathcal{H} \to \mathbb{B}^n$ be a conformal homeomorphism such that $\phi(0)=0$. Let $\ep>0$. Then $\phi^{-1}(\ov B(0;1-\ep))$ is compactly contained in $\mathcal{H}$ and hence is contained in $\Omega_j$ for all large $j$. Since $\Omega_j \subset \mathcal{H}$, $\phi$ is a $1$-quasiconformal map of $\Om_j$ into $\mathbb{B}^n$. Thus $B(0;1-\ep) \subset \phi(\Om_j)$ and hence the map $\phi$ is a candidate for $s^1_{\Omega_j}(0)$ for all large $j$. Therefore,
\[
1-\ep \leq s^1_{\Omega_j}(0),
\]
for all large $j$. It follows that
\[
1-\ep \leq \liminf_{j \to \infty}s^1_{\Omega_j}(0),
\]
and hence
\[
1\leq \liminf_{j \to \infty}s^1_{\Omega_j}(0),
\]
as required.






\end{proof}
\subsection{Some examples of $K$-squeezing functions} It has been previously seen that $s^K_{\mbb B^n} \equiv 1$ for all $k \ge 1$. We present two examples to illustrate the behaviour of $s^K_{\Om}$ when $\Om$ is not conformally equivalent to $\mbb B^n$.

\begin{prop}
Consider the ellipsoid
\[
\mathcal E=\left\{x \in \mathbb{R}^n: \sum_{j=1}^{n}\frac{x_j^2}{a_j^2}<1\right\},
\]
where $n\ge 3$, $0<a_j<\infty$, $1 \leq j\leq n$, and suppose that $a_i \neq a_j$ for some pair $i \neq j$. Then there exists $K_0 \in (1,\infty)$ such that
\begin{itemize}
\item [(a)] for every $K \in [1,K_0)$, $s^K_{\mathcal E}<1$.
\item [(b)] for every $K\in [K_0,\infty)$,  $s^K_{\mathcal E}\equiv 1$.
\end{itemize}
\end{prop}
\begin{proof}
Since $\mathcal E$ is the image of $\mathbb{B}^n$ under an invertible linear map, it is quasiconformally equivalent to $\mathbb{B}^n$. Let
\[
K_0=\inf\big\{K \in [1,\infty): \text{$\mathcal E$ is $K$-quasiconformally equivalent to $\mathbb{B}^n$}\big\}.
\]
Then there exists a $K_0$-quasiconformal map of $\mathcal E$ onto $\mathbb{B}^n$ (see, for instance, \cite[Theorem~7.1.1]{GMP}). Since conformal maps in dimension $\geq 3$ are M\"obius transformations that map balls to balls or half-spaces, we must have $K_0>1$. Thus (a) and (b) follow from Proposition~\ref{s-ball} and Proposition~\ref{s-dec}.
\end{proof}

Thus, for the ellipsoids $\mathcal E$, the family of conformal invariants $\{s^K_{\mathcal E}: K \ge 1\}$ undergoes a transition at an intermediate critical dilatation $K_0 > 1$. For dilatations greater than the critical value $K_0$, the $K$-squeezing functions are all identically constant and take the maximal value $1$, whereas for lower dilatations, they are strictly less than $1$. On the other hand, for each $1 \le K < K_0$, $s^K_{\mathcal E} \ra 1$ near the boundary $\pa \mathcal E$. Thus the invariants $s^K_{\mathcal E}$ are non-constant for $1 \le K < K_0$. 

\medskip

Recall from Section~2 that $\phi_{K,n}:(0,1) \to (0,1)$ is the function defined by 
\[
\phi_{K,n}(r)=\frac{1}{\ga_n^{-1}\left(K\gamma_n\left(\frac{1}{r}\right)\right)}.
\]
\begin{prop}
Let $\Omega=\mathbb{B}^n \setminus \{0\}$ and $K \in [1,\infty)$. Then
\begin{equation}\label{s_pb-bd}
\vert x \vert \leq s^K_{\Omega}(x) \leq \phi_{K,n}\big(\vert x \vert\big), \quad x \in \Omega.
\end{equation}
In particular, $s^1_{\Omega}(x)=\vert x\vert$ for all $x \in \Omega$.
\end{prop}
\begin{proof}
Fix $K \in [1, \infty)$ and $x \in \Omega$. Let
\[
T_x=p_x \circ \sigma_x
\]
where $\sigma_x$ is the inversion in the Euclidean sphere of radius $\sqrt{\vert x \vert^{-2}-1}$ centred at $x/\vert x\vert^2$ and $p_x$ is the reflection in the hyperplane through the origin orthogonal to $x$. Then $T_x$ is a sense-preserving M\"obius transformation of $\overline{\mathbb{R}^n}$ that preserves $\mathbb{B}^n$ (see \cite[1.34]{V}). Moreover, $T_x(x)=0$ and $T_x(0)=-x$. Letting $f=T_x\vert_{\Omega}$, we see that $f: \Omega \to \mathbb{B}^n$ is a sense-preserving conformal map such that $f(x)=0$ and
\[
f(\Omega)=\mathbb{B}^n\setminus \{-x\} \supset  B\big(0;\vert x\vert \big).
\]
Therefore
\begin{equation}\label{s_pb-lb}
s^K_{\Omega}(x) \geq s^1_{\Omega}(x) \geq \vert x\vert.
\end{equation}

On the other hand, let $g:\Omega \to \mathbb{B}^n$ be a $K$-quasiconformal map such that $g(x)=0$ and $B(0;s^K_{\Omega}(x)) \subset g(\Omega)$. Then, by \cite[Theorem~17.3]{V}, $g$ has a limit $y$ at $0$ and $y$ is an isolated boundary point of $g(\Omega)$. Moreover, the map $G$ defined by $G\vert_{\Omega}=g$ and $G(0)=y$ is a $K$-quasiconformal map of $\mathbb{B}^n$ onto $g(\Omega) \cup \{y\}$. Since $g(\Omega) \cup \{y\}$ must be an open subset of $\overline{\mathbb{R}^n}$, we have $g(\Omega) \cup \{y\} \subset \mathbb{B}^n$. Now
\[
\mu_{\mathbb{B}^n}\big(G(x),G(0)\big) \leq  \mu_{G(\mathbb{B}^n)}\big(G(x),G(0)\big) \leq K \mu_{\mathbb{B}^n}(x,0)
\]
which implies that
\[
\ga_n\left(\frac{1}{\vert y\vert}\right) \leq K \ga_n\left(\frac{1}{\vert x\vert}\right),
\]
and hence $\vert y \vert \leq \phi_{K,n}(\vert x\vert)$. As $y$ is a boundary point of $g(\Omega)$ and $B(0;s^K_{\Omega}(x)) \subset g(\Omega)$, it follows that
\begin{equation}\label{s_pb-ub}
s^K_{\Omega}(x) \leq \vert y\vert \leq \phi_{K,n}\big(\vert x \vert\big).
\end{equation}
Combining \eqref{s_pb-lb} and \eqref{s_pb-ub}, we obtain \eqref{s_pb-bd}.
\end{proof}

The following corollary is immediate from \eqref{phi-limit}.

\begin{cor}
Let $\Omega=\mathbb{B}^n \setminus \{0\}$ and $K \in [1,\infty)$. Then
\[
\lim_{x \to 0}s^K_{\Omega}(x)=0.
\]
\end{cor}

More generally, we have the following estimate:

\begin{prop}\label{null-set-sq}
Let $\Omega \subset \mathbb{R}^n$ be a bounded domain and $\Omega'=\Omega \setminus P$, where $P$ is a compact subset of $\Omega$ of capacity zero. Then for $K \in [1,\infty)$,
\[
s^K_{\Omega'}(x) \leq \frac{1}{\ga_n^{-1}\big(K\mu_{\Omega}(x,P)\big)}, \quad x \in \Omega'.
\]
\end{prop}
\begin{proof}
Fix $x \in \Omega'$. Let $f: \Omega' \to \mathbb{B}^n$ be a $K$-quasiconformal map such that $f(0)=x$ and $B(0;r) \subset f(\Omega')$. We have to prove that
\[
r \leq \frac{1}{\ga_n^{-1}\big(K\mu_{\Omega}(x,P)\big)}.
\]
This will follow at once if we show that $d(0,\pa f(\Omega'))$ is bounded above by the right-hand side. By Theorem~6.4.26 of \cite{GMP}, $f$ extends to a $K$-quasiconformal map $F$ of $\Omega$ onto $f(\Omega') \cup C_{f}(P)$, where $C_{f}(P)$ is the cluster set of $P$ under $f$. Evidently, $C_{f}(P)\subset \overline{\mathbb{B}^n}$ and thus $F: \Omega \to \overline{\mathbb{B}^n}$. Since $F(\Omega)$ is an open subset of $\overline{\mathbb{R}^n}$, we must have $F: \Omega \to \mathbb{B}^n$. Also, note that $F(P) \subset \partial f(\Omega')$ and $F(P)$ is a compact subset of $\mathbb{B}^n$. Therefore there exists $q \in F(P)$ such that
\[
\mu_{\mathbb{B}^n}(0,q)=\mu_{\mathbb{B}^n}\big(0,F(P)\big)=\mu_{\mathbb{B}^n}\big(F(x),F(P)\big)\leq \mu_{F(\Omega)}\big(F(x),F(P)\big) \leq K \mu_{\Omega}(x,P).
\]
Since $\mu_{\mathbb{B}^n}(0,q)=\ga_n(1/\vert q\vert)$ and $\ga_n$ is strictly decreasing, we obtain
\[
\vert q \vert \leq \frac{1}{\ga_n^{-1}\big(K \mu_{\Omega}(x,P)\big)}.
\]
As $q \in F(P) \subset \partial f(\Omega')$, the above implies that $d(0,\partial f(\Omega'))$ is bounded above by the right-hand side as required.
\end{proof}

As we shall see later, the above result can also be recovered as a corollary of Propositions \ref{null-set-frid} and \ref{sq-frid}. The following corollary is immediate from \eqref{ga-dec}:

\begin{cor}
Under the hypothesis of Proposition~\ref{null-set-sq},
\[
\lim_{x \to P} s^K_{\Omega'}(x)=0.
\]
\end{cor}

\section{The $K$-Fridman invariant}
Let $M$ be a $C^1$-smooth $n$-dimensional Riemannian manifold of Class II. 
Fix $ x \in M  $ and $ K > 1 $. Choose $ \epsilon > 0 $ such that $ (1 + \epsilon)^{2(n-1)} < K $. Since $ M $ is $ C^1$-smooth, by \cite[Lemma 3.1(a)]{F4}, there exist a neighbourhood $ U $ of $ x $ and a $ (1 + \epsilon)$-bilipschitz homeomorphism $ \theta $ of $ U $ onto a domain $ U' $ in $ \widehat{\mathbb{R}}^n $. As a consequence, both $ \theta, \theta^{-1} $ are $ (1 + \epsilon)^{2(n-1)} $-quasiconformal. Since $ U' $ is a domain containing $ \theta(x) $, choose a Euclidean ball $ B \big(\theta(x); s \big) $ whose closure is contained in $ U' $. Let $ T:\mathbb{B}^n \to B \big(\theta(x); s \big) $ be the affine similarity defined by 
\[ T( y )= \theta(x) + s y.
\]
Then $ f:= \theta^{-1} \circ T : \mathbb{B}^n \to M $ is a $ K$-quasiconformal map with $ f(0) = x $. Moreover, $ 
f(\mathbb{B}^n) = \theta^{-1} \big( B \big(\theta(x); s \big) \big)
$ is a neighbourhood of $ x $. Since $ M $ is of Class II, $ \mu_M $ induces the manifold topology, and hence there exists a $ r > 0 $ such that 
\[
B_M (x;r) \subset f(\mathbb{B}^n), 
\]
where $ B_M (x; r) $ is the open ball in $M$ of radius $ r >0 $, centered at $ x $, with respect to the metric $ \mu_M $.

\medskip

\noindent Define 
\begin{equation}\label{h^K}
h^K_M(x) = \inf \{ 1/r: B_M(x;r) \subset f(\mbb B^n): f \in QC_K(\mbb B^n, M), f(0) = x \} \, .
\end{equation}
The definition of $ h_M^K$ is motivated by the analogous function introduced in \cite{Fridman1}, \cite{Fridman2} for Kobayashi hyperbolic manifolds in the holomorphic category, with respect to the Kobayashi metric. The function defined in (\ref{h^K}) will henceforth be referred to as the $ K $-Fridman function.

\medskip

By the preceding argument, the family of maps that define this invariant is non-empty for each $ K > 1 $. However, in general, it may be empty when $K = 1$. In contrast, when $ M $ is a bounded domain in $ \mathbb R^n $, it is evident that this family of maps is non-empty for every $ K\geq 1 $.

\medskip 

Let $M'$ be another $C^1$-smooth $n$-dimensional Riemannian manifold of Class II and $\phi: M \to M'$ an $L$-quasiconformal homeomorphism. Suppose $1/r$ is a candidate for $h^K_M(x)$. Then there exists a $K$-quasiconformal map $f: \mathbb{B}^n \to M$ such that $f(0)=x$ and $B_M(x;r) \subset f(\mathbb{B}^n)$. Thus, $\phi \circ f: \mathbb{B}^n \to M'$ is a $KL$-quasiconformal map such that $\phi\circ f(0)=\phi(x)$ and
\[
B_{M'}(\phi(x); r/L) \subset \phi(B_M(x;r)) \subset \phi \circ f(\mathbb{B}^n).
\]
This shows that $ L/r$ is a candidate for $h^{KL}_{M'}(\phi(x))$. Therefore
\[
h^{KL}_{M'}\big(\phi(x)\big) \leq {L}\; h^K_M(x).
\]
In particular, if $\phi$ is conformal, then
\[
h^K_{M}(x)=h^K_{M'}\big(\phi(x)\big).
\]
Moreover, it follows directly from the definition that for $ x \in M $
\[
h^{K_2}_M(x) \leq h^{K_1}_M(x),
\]
whenever $K_1 \leq K_2$.

\subsection{Examples} The following two examples help to set the tone for what can be expected for the $K$-Fridman invariant.

\medskip

\begin{prop}
For the unit ball $\mathbb{B}^n \subset \mbb R^n$, $h^K_{\mathbb{B}^n}\equiv 0$.
\end{prop}
\begin{proof}
Indeed, it suffices to show that $h^K_{\mathbb{B}^n}(0)=0$, as the group of M\"obius transformations acts transitively on $\mathbb{B}^n$, and these transformations are conformal. To establish this, note that from \eqref{mu-ball}-\eqref{ga-dec}, given any $r>0$, the following inequalities are equivalent:
\[
\mu_{\mathbb{B}^n}(y,0)< r \ , \ \ga\left(\frac{1}{\vert y\vert }\right)<r \ , \ \vert y \vert < \frac{1}{\gamma^{-1}(r)}. 
\]
Since $ \gamma^{-1}:(0, \infty) \to (1, \infty) $,
it follows that
\[
B_{\mathbb{B}^n}(0;r) \subset \mathbb{B}^n.
\]
As the identity map $\mathbf{1}_{\mathbb{B}^n}$ is $K$-quasiconformal for every $K \ge 1$, it follows that $1/r$ is a candidate for $h^K_{\mathbb{B}^n}(0)$ for every $r>0$ and hence
\[
h^K_{\mathbb{B}^n}(0)=0
\]
as required.
\end{proof}

\medskip

\begin{prop}
For the punctured unit ball $ \mathbb{B}^n \setminus \{0\} $, 
\[
h^K_{\mathbb{B}^n \setminus \{0\}} (x) = \frac{1}{ \gamma_n \left( \frac{1}{|x|} \right)}, 
\]
for $ x \in \mathbb{B}^n \setminus \{0\} $. In particular, $ h^K_{\mathbb{B}^n \setminus \{0\}} (x) \to \infty $ as $ x \to 0 $ and $ h^K_{\mathbb{B}^n \setminus \{0\}} (x) \to 0 $ as $ x \to \partial \mathbb{B}^n $. 
\end{prop}

\begin{proof}
To see this, fix $ \epsilon > 0 $ and let  $f : \mathbb{B}^n \rightarrow \mathbb{B}^n \setminus \{0\} $ be a $K$-quasiconformal embedding with $ f(0) = x $ and 
$ B_{\mathbb{B}^n \setminus \{0\}} (x, r) \subset f( \mathbb{B}^n ) $ such that 
\[
\frac{1}{r} < h^K_{\mathbb{B}^n \setminus \{0\} } (x) + \epsilon.
\]
Then the claim is that $ r \leq \gamma_n \left( \frac{1}{|x|} \right) $. Assume, on the contrary, that 
\[
 r > \gamma_n \left( \frac{1}{|x|} \right) = \mu_{\mathbb{B}^n} (x,0).
\]
Since the function $ y \mapsto \mu_{\mathbb{B}^n} (x,y) $ is continuous at the origin, there exists a $ \delta > 0 $ such that $ \mu_{\mathbb{B}^n} (x,y)
< r $ whenever $ |y | < \delta $. Hence
\[
B (0;\delta ) \subset  B_{\mathbb{B}^n} (x; r). 
\]
Now, using the fact that $ \{0\} $ has zero $n$-capacity, we have 
\[ 
\mu_{\mathbb{B}^n} (s,t) = \mu_{\mathbb{B}^n \setminus \{0\} } (s,t) 
\]
for all $ s, t \in \mathbb{B}^n \setminus \{0\} $. As a consequence,
\[
 B_{\mathbb{B}^n \setminus \{0\}} (x; r) =  B_{\mathbb{B}^n } (x; r) \setminus \{0\},
\]
which, in turn, implies that 
\[
B (0;\delta )  \setminus \{0\} \subset  B_{\mathbb{B}^n} (x; r) \setminus \{0\} = B_{\mathbb{B}^n \setminus \{0\}} (x; r) \subset f(\mathbb{B}^n) . 
\]
The inclusion $ B (0;\delta ) \setminus \{0\} \subset  f(\mathbb{B}^n) $ allows us to consider the restriction of the inverse map $ f^{-1} $ to $ B (0;\delta ) \setminus \{0\} $. Moreover,
\[
g:= f^{-1}: B (0;\delta ) \setminus \{0\} \to \mathbb{B}^n
\]
is a $ K$-quasiconformal mapping. In this setting, Corollary 6.4.27 of \cite{GMP} guarantees the existence of the limit 
\[
\lim_{y \to 0} g (y) = x_0,
\]
where $ x_0 $ is an isolated boundary point of $  g \big( B (0;\delta ) \setminus \{0\} \big) $. Moreover, the extension $ G $
of $ g $ to $ B (0;\delta )  $ defined at the origin by $ G ( 0) = x_0 $ is a $ K $-quasiconformal homeomorphism from  $ B (0;\delta )  $ onto the domain 
\[
X= g \big( B (0;\delta ) \setminus \{0\} \big) \cup \{x_0\} \subset \overline{\mathbb{B}^n}.
\]
Since $ X $ is an open subset of $ \mathbb{R}^n $ contained in $ \overline{\mathbb{B}^n} $, every point of $ X $ lies in $ \mathbb{B}^n $. In particular, $ x_0 \in  \mathbb{B}^n $.

\medskip

Next, let $ \{ y_n \} $ be a sequence in $B (0;\delta ) \setminus \{0\} $ such that $ y_n \to 0 $. Since $ G $ is continuous at the origin and $ G(0)= x_0 $, we have $ G(y_n) \to x_0 $. As $ f $ is continuous at $ x_0 $, it follows that
$ (f \circ G)(y_n) \to f(x_0) $. On the other hand, 
\[
( f \circ G)(y_n)= (f \circ g) (y_n) = y_n \to 0.
\]
Therefore $ f(x_0)=0 $. This contradicts the fact that $ f(x_0)\in \mathbb B^n\setminus \{0\} $. Hence, the assumption $ r > \gamma_n \left( \frac{1}{|x|} \right) $ is false, and therefore  $ r \leq \gamma_n \left( \frac{1}{|x|} \right) $. As a consequence,
\[
h^K_{\mathbb{B}^n \setminus \{0\} } (x) + \epsilon > \frac{1}{r} \geq \frac{1}{\gamma_n \left( \frac{1}{|x|} \right)}.
\]
Since $ \epsilon > 0 $ was arbitrary, it follows that $ h^K_{\mathbb{B}^n \setminus \{0\} } (x) \geq \frac{1}{\gamma_n \left( \frac{1}{|x|} \right)} $.

\medskip

To establish the upper bound $ h^K_{\mathbb{B}^n \setminus \{0\} } (x) \leq \frac{1}{\gamma_n \left( \frac{1}{|x|} \right)} $, recall that
\[ 
\mu_{\mathbb{B}^n} (s,t) = \mu_{\mathbb{B}^n \setminus \{0\} } (s,t) 
\]
for all $ s, t \in \mathbb{B}^n \setminus \{0\} $.
Setting $ R := \gamma_n \left( \frac{1}{|x|} \right) = \mu_{\mathbb{B}^n} (x,0) $, it is evident that $ 0 \notin B_{\mathbb{B}^n} (x;R) $,
and therefore
\[
 B_{\mathbb{B}^n \setminus \{0\}} (x; R) =  B_{\mathbb{B}^n } (x; R).
 \]
 Now, by definition, 
\begin{align*}
B_{\mathbb{B}^n} (x; R) = & \left\{ y \in \mathbb{B}^n: \gamma_n \left(\frac{1}{\tanh \frac{1}{2} \rho_{\mathbb{B}^n}(y,x)}\right) < R \right\},
\\
= & \left\{ y \in \mathbb{B}^n: { \rho_{\mathbb{B}^n}(y,x)} < 2 \tanh^{-1} \left(|x| \right) \right\}, \\
= & \left\{ y \in \mathbb{B}^n: { \rho_{\mathbb{B}^n}(y,x)} < \rho_{\mathbb{B}^n}(0,x) \right\},
\end{align*}
that is, $ B_{\mathbb{B}^n} (x, R) $ is the $ \rho_{\mathbb{B}^n}$-ball in $ \mathbb{B}^n $ with center $ x $ and radius $ \rho_{\mathbb{B}^n}(0,x) $. It follows from the explicit expression of the hyperbolic distance on $ \mathbb{B}^n $ that $ B_{\mathbb{B}^n} (x; R) $ is a Euclidean ball contained in $ \mathbb{B}^n \setminus \{0\}$.
As a consequence, there is a M\"{o}bius transformation 
\[
f: \mathbb{B}^n \rightarrow  B_{\mathbb{B}^n} (x; R)  
\]
from $ \mathbb{B}^n $ onto $ B_{\mathbb{B}^n} (x; R)  $ with $ f(0) = x $. Since $ f $ is conformal, it is $ 1$-quasiconformal and hence $ K $-quasiconformal for every $ K \geq 1 $. Hence $ f $ is a candidate map in the definition of $ h^K_{\mathbb{B}^n \setminus \{0\} } (x) $ with
\[
B_{\mathbb{B}^n \setminus \{0\}} (x; R) = f(\mathbb{B}^n) \subset \mathbb{B}^n \setminus \{0\},
\]
which, in turn, implies that
\[
 h^K_{\mathbb{B}^n \setminus \{0\} } (x) \leq \frac{1}{R} = \frac{1}{\gamma_n \left( \frac{1}{|x|} \right)}.  
\]
and this completes the claim.
\end{proof}

The above example shows that, unlike the $K$-squeezing function, the $K$-Fridman invariant $h^K_{M}$ is, in general, not strictly increasing in the parameter $K$. It also motivates the following estimate for the $K$-Fridman invariant on complements of capacity zero sets.

\begin{prop}\label{null-set-frid}
Let $\Omega \subset \mathbb{R}^n$ be a bounded domain and $\Omega'=\Omega \setminus P$, where $P \subset \Omega$ is a compact set of zero capacity. Then for $K \ge 1$,
\[
h^K_{\Omega'}(x) \geq \frac{1}{\mu_{\Omega}(x,P)}, \quad x \in \Omega'.
\]
\end{prop}

\begin{proof}
Fix $x \in \Omega'$. Let $f: \mathbb{B}^n \to \Omega'$ be a $K$-quasiconformal embedding such that $f(0)=x$ and $B_{\Omega'}(x;r) \subset f(\mathbb{B}^n)$. We have to prove that $r \leq \mu_{\Omega}(x,P)$.

\medskip

If possible, assume that $r>\mu_{\Om}(x,P)$. Since $P$ is a compact subset of $\Omega$, there exists $p \in P$ such that $\mu_{\Omega}(x,P)=\mu_{\Omega}(x,p)$. Note that $p \in \partial \Omega'$. By continuity of $\mu_{\Omega}(x,\cdot)$, there exists $\de>0$ such that for all $y\in B(p;\de)$, $\mu_{\Omega}(x,y)<r$. Therefore
\[
B(p;\delta) \subset B_{\Omega}(x;r).
\]
As $\mu_{\Omega'}=\mu_{\Omega}\vert_{\Omega'\times \Omega'}$,
\[
B_{\Omega'}(x;r)=B_{\Omega}(x;r)\setminus P,
\]
and hence
\[
B(p;\delta) \setminus P \subset B_{\Omega'}(x;r) \subset f(\mathbb{B}^n).
\]
Let $g=f^{-1}\vert_{B(p;\delta) \setminus P }$. Then $g$ extends to a $K$-quasiconformal map $G$ of $B(p;\delta)$ into $\overline{\mathbb{B}^n}$ (see Theorem~6.4.26 of \cite{GMP}). Since $G(B(p;\delta))$ is an open set in $\overline{\mathbb{R}^n}$, we must have $G(B(p;\delta)) \subset \mathbb{B}^n$. In particular, $G(p) \in \mathbb{B}^n$. This implies $f\circ G(p) \in \Omega'$.

\medskip

On the other hand, since $p \in \partial \Omega'$, there exists a sequence $ \{y_n \}$ in $\Omega'$ that converges to $p$. Thus, $g(y_n)=G(y_n) \to G(p)$ implies  $f\circ g(y_n) \to f\circ G(p)$, which in turn implies $y_n \to f\circ G(p)$ and, finally, $f\circ G(p)=p$. This is a contradiction since $f \circ G(p) \in \Om'$ whereas $p \in \partial \Omega'$.
\end{proof}

\medskip
\subsection{Existence of extremals and continuity of $h^K_{M}$}

We first state a version of Lemma 1.1 of \cite{Fridman2} adapted to our setting.

\begin{prop} \label{P1}
Let $M$ be a $C^1$-smooth $n$-dimensional Riemannian manifold of Class II, and let $ D $ be a bounded domain in $ \mathbb{R}^n$. Suppose that there exist two relatively compact sets $ K_1 \subset D $ and $ K_2 \subset M $ and a sequence $ \{F_k \} $ of mappings $ F_k : D \rightarrow F_k (D) \subset M $
satisfying the following conditions:
 
\begin{enumerate}
    \item [(i)] there exists a $ K \geq 1 $, independent of $ k $, such that for every $ k \geq 1 $, the mapping $ F_k : D \rightarrow F_k (D) $ is a $K$-quasiconformal homeomorphism,
     \item [(ii)] for each $ k \geq 1 $, there exists a point $ z_k \in K_1 $ such that $ F_k (z_k ) \in K_2 $,
     \item [(iii)] for any compact $ L \subset M $, there exists a number $ s = s(L) $ such that $ L \subset F_s (D) $.
\end{enumerate}
 
Then $ M $ is K-quasiconformally equivalent to $D$.    
\end{prop}

The above proposition follows from the arguments presented in the proof of Theorem 3.2(i) in \cite{BMV}. We therefore do not repeat the proof here.

\begin{prop}
Let $M$ be a $C^1$-smooth $n$-dimensional Riemannian manifold of Class II and $ K \in [1, \infty)$. 

\begin{enumerate}
\item [(i)] If there is a $x_0 \in M $ such that $h^K_{M}(x_0)=0$, then $ M $ is $K$-quasiconformally equivalent to the unit ball of $ \mathbb{R}^n $ and $h^K_{M} \equiv 0$.

\item [(ii)] If there is a $x_0 \in M $ such that $h^K_{M}(x_0) > 0$ and $ (M, \mu_M) $ is proper, then there is a $K-$quasiconformal mapping $ f_0: \mathbb{B}^n \to M $ with $ f_0(\mathbb{B}^n) \supset B_M \big( x_0, 1/h^K_M(x_0) \big) $.

\item [(iii)] We have that $h^K_M (\cdot) $ is continuous on $ M $. 
\end{enumerate}

\end{prop}

\begin{proof} We first prove $(i)$. By the definition of $ h^K_M (\cdot) $,
there is a sequence $ \{ \phi_k \} $ of $ K$-quasiconformal mappings $ \phi_k : \mathbb{B}^n \rightarrow M $ such that $ \phi_k(0) = x_0 $ and $ B_M( x_0; k) \subset \phi_k(\mathbb{B}^n ) $.

\medskip

Then, for every compact $ L\subset M $, the continuity of $ \mu_M(x_0,\cdot) $ gives $ c=c(L)>0 $ such that $ L \subset B_M(x_0,c) $. It follows that  $ L \subset \phi_k(\mathbb B^n) $ for every $ k \geq c $, so the domains $  \phi_k(\mathbb B^n) $ exhaust $ M $.

\medskip

Applying Proposition \ref{P1} to the sequence $ \{ \phi_k \} $ with $ D = \mathbb{B}^n $, $ K_1 = \{ 0 \} $ and $ K_2 = \{ x_0 \} $, it follows that 
$ M $ is $K$-quasiconformally equivalent to $ \mathbb{B}^n $.  
\medskip

Let $f: \mathbb{B}^n \to M $ be a $K$-quasiconformal homeomorphism from $ \mathbb{B}^n $ onto $ M $. Fix $ y \in M $, and let $ x \in \mathbb B^n $ satisfy $ f(x)=y $. Let $ T $ be a conformal automorphism of $ \mathbb B^n $ such that $ T(0)=x $. Then $ f \circ T: \mathbb B^n \to M $ is a $ K$-quasiconformal homeomorphism with $ (f\circ T)(0)=y $. Since
\[ 
B_M(y;r) \subset M =(f\circ T)(\mathbb B^n),
\]
for every $ r>0 $, it follows that $ 1/r $ is admissible in the definition of $ h_M^K(y)$. Hence $ h_M^K(y)=0 $. Since $ y \in M $ was arbitrary, $ h^K_{M} \equiv 0 $.

\medskip

For $(ii)$, by definition, there exist sequences of $ K$-quasiconformal embeddings $ f_k: \mathbb{B}^n \to M $ and positive real numbers $ r_k > 0 $ with $ f_k(0) = x_0 $, $ f_k( \mathbb{B}^n) \supset B_M(x;r_k) $ and $ 1/r_k \to h_M^K(x_0) $ as $ k \to \infty $. 

\medskip

It follows from Theorem C(b) of \cite{F1} that $ \{ f_k\} $ is $K$-Lipschitzian in the metrics $ \mu_{\mathbb{B}^n} $, $ \mu_{M} $
and hence the family $ \{ f_k\} $ is uniformly equicontinuous on $ \mathbb{B}^n $.
The claim is that $ \{f_k\} $ is normal. To establish this, observe that $ \{f_k\} $ is locally uniformly bounded on $ \mathbb{B}^n $. Indeed, for a fixed compact subset $ L $ of $ \mathbb{B}^n $, 
\[
\mu_M\big( f_k(z), x_0 \big) = \mu_M\big( f_k(z), f_k(0) \big) \leq K \mu_{\mathbb{B}^n} ( z, 0 )
\]
for every $ z \in L $. Since $ L $ is compact, there exists a constant $ C > 0 $, independent of $ k $, such that $ \mu_{\mathbb B^n}(z,0)\le C $ for all  $ z\in L $. As a consequence, 
\[
f_k(L) \subset \overline{B}_M (x_0; KC) 
\]
for every $ k $. The properness of $ (M, \mu_M) $ ensures that the closed ball $ \overline{B}_M (x_0, KC) $ is compact, which, in turn, implies that $ \{ f_k \} $ is uniformly bounded on $ L $. An application of the Arzel\`{a}–Ascoli theorem yields a subsequence of $ \{ f_k \} $ that converges uniformly on $ L $. A diagonal argument then gives a subsequence that converges locally uniformly on $ \mathbb B^n $, establishing the normality of the family $ \{f_k \}$.

\medskip

Passing to a subsequence, we may assume that $ \{f_k \} $ converges uniformly on compact subsets of $ \mathbb{B}^n $ to a continuous map $ f_0 $. By Theorem 7.3 of \cite{F1}, the limit $f_0 $ is a constant or a $K$-quasiconformal mapping from $ \mathbb{B}^n $ to $ M $. Since $ f_0(0) = x_0 $, it follows that $ f_0(\mathbb{B}^n) \subset M $ in either case. We shall show that $f_0$ realizes $h_M^K(x_0)$.

\medskip

To this end, fix $ \epsilon > 0 $ sufficiently small and consider the $ K$-quasiconformal mappings
\[
g_k := (f_k)^{-1}: B_M \left(x_0; \frac{1}{h_M^k(x_0)} - \epsilon \right) \to \mathbb{B}^n.  
\]
It follows from Theorem 3.1 of \cite{BMV} that $ \{g_k\} $ is a normal family. Let $ g_0 $ be a subsequential limit of $\{g_k\} $. Since $ g_k(x_0) = 0 $ for every $ k $, by the same reasoning as above, it follows that 
\[
g_0: B_M \left(x_0; \frac{1}{h_M^k(x_0)} - \epsilon \right) \to \mathbb{B}^n. 
\]
Now, for any $ x \in B_M \left(x_0; \frac{1}{h_M^k(x_0)} - \epsilon \right) $, there is a compact set $ L' \Subset \mathbb{B}^n $ such that $ g_k(p) \in L' $ for all $ k $ large. Therefore
\[
x = f_k \circ g_k (x) \subset f_k(L').
\]
Let $ y_k \in L' $ be such that $ f_k(y_k) = x $. After passing to a subsequence, we may assume that $ y_k \to y_0 \in L'$. Since $ \{ f_k \} $ converges to $ f_0 $ uniformly on compact subsets of $ \mathbb{B}^n $, it follows that
\[
 f_k(y_k) \to f_0(y_0) 
\]
as $ k \to \infty $. Therefore $ x = f_0(y_0) \in \mathbb{B}^n $. This is exactly the assertion that
\[
B_M \left(x_0; \frac{1}{h_M^k(x_0)} - \epsilon \right)  \subset f_0(\mathbb{B}^n). 
\]
Since this inclusion holds for every sufficiently small $ \epsilon >0 $, letting $ \epsilon \ra 0 $ yields
\[
B_M \left(x_0; \frac{1}{ h_M^K(x_0) } \right)
\subseteq f_0(\mathbb B^n).
\]
Since $ f_0( \mathbb{B}^n) $ contains a non-empty open set, $ f_0 $ must be non-constant. It follows from the preceding dichotomy that $ f_0 : \mathbb{B}^n \to M $ is $ K$-quasiconformal. This establishes that $ f_0 $ is an extremal map for the variational problem defining $ h_M^K(x_0) $.

\medskip

Part $(iii)$ follows by making the relevant changes in the proof of Theorem 1.3 of \cite{Fridman2} using the continuity of $ \mu_M $ and Theorem $7.5$ of \cite{F5}.
\end{proof}

\subsection{Dependence of $h^K_{\Om}$ on $\Om$}

We investigate the relationship between the $K$-Fridman functions of a monotone sequence of domains and that of its limit domain. We begin with the case of an increasing sequence of domains.

\begin{thm} \label{inc}
Let $ \Omega \subset \mathbb{R}^n$ be a bounded domain and $ \{ \Omega_j \}$ be an increasing sequence of sub-domains of $ \Omega $ such that $\cup_{j} \Omega_j= \Omega $. Then, for $ K \ge 1$, 
\[
\lim_{j \to \infty} h^K_{\Omega_j}(x) = h^K_{\Omega}(x),
\]
for every $ x \in \Omega $. 
\end{thm}

\begin{proof} 
Let $ r > 0 $ be such that $ h^K_{\Omega}(x_0) < 1/r $. Then there is a $ K$-quasiconformal embedding $ f : \mathbb{B}^n \to \Omega $ such that $ f(0) = x_0 $ and $ f(\mathbb{B}^n) \supset B_\Omega(x_0; r) $. 

\medskip

For $ \epsilon > 0 $, there exists a $ \delta > 0 $ such that $ f \big( B(0; 1- \delta) \big) \supset B_\Omega(x_0; r - \epsilon) $. Since $ f \big(  B(0; 1- \delta)\big) $ is relatively compact in $ \Omega $, $ f \big(  B(0; 1- \delta) \big) \subset \Omega_j $ for all $ j $ sufficiently large. The inclusion $ \Omega_j \subset \Omega $ guarantees that $ B_{\Omega_j}(x_0; r - \epsilon) \subset B_\Omega(x_0; r ) $ for all $ j $ sufficiently large. It follows that $ g: \mathbb{B}^n \to \Omega_j $ defined by setting $ g(x) = f \big( (1- \delta) x \big) $ is a $ K$-quasiconformal mapping with $ g(0) = x_0 $ and 
\[
g( \mathbb{B}^n) = f \big(  B(0; 1- \delta) \big) \supset B_\Omega(x_0; r - \epsilon) \supset B_{\Omega_j}(x_0; r - \epsilon),
\]
for all $ j $ sufficiently large. Therefore, $ h^K_{\Omega_j}(x_0) \leq 1/(r - \epsilon) $ for all $ j $ sufficiently large. Letting $ \epsilon \to 0 $, it follows that $ \limsup_{j \to \infty} h^K_{\Omega_j}(x_0) \leq 1/r $. Since the above estimate holds for every $ r>0 $ such that $ h^K_\Omega(x_0) < 1/r $, letting $ 1/r $ decrease to $ h^K_\Omega(x_0) $ yields
\[
\limsup_{j \to \infty} h^K_{\Omega_j}(x_0) \leq h^K_{\Omega}(x_0).
\]
In order to establish $ \liminf_{j \to \infty} h^K_{\Omega_j}(x_0) \geq h^K_{\Omega}(x_0) $, the following observation will be needed:

\medskip

\noindent {\it Claim}: For any $ R > 0 $ and $ \epsilon > 0 $ sufficiently small,
\begin{equation} \label{inclu}
 B_\Omega(x_0; R - \epsilon) \subset B_{\Omega_j} (x_0; R)   
\end{equation}
for all $ j $ sufficiently large. 

\medskip

To verify this, fix $ R > 0 $ and let $ x \in \Omega $ be such that $ \mu_\Omega(x, x_0) < R - \epsilon $. Then there is a compact continuum  $ C \subset \Omega $ that contains $ x $ and $ x_0 $ such that 
\begin{equation} \label{ineq}
\mbox{Cap}_\Omega(C) < \mu_\Omega(x, x_0) + \epsilon /2 < R - \epsilon/2.
\end{equation}
Since $ \Omega = \cup_{j} \Omega_j $, $ C \Subset \Omega_j $ for all $ j $ sufficiently large. Moreover, since the domains $ \{ \Omega_j \} $ form an increasing sequence, the sequence $ \{ \mbox{Cap}_{\Omega_j}(C) \} $ is decreasing and bounded below by $ \mbox{Cap}_{\Omega}(C) $. Setting $ \ell = \lim_{j \to \infty} \mbox{Cap}_{\Omega_j}(C) $, it is evident that $ \ell \geq \mbox{Cap}_{\Omega}(C) $. 

\medskip

We show that $ \ell = \mbox{Cap}_{\Omega}(C) $. Indeed, let $ u \in H_0(\Omega) $ be admissible for $ C $. By the definition of $ H_0(\Omega) $, $ u $ has compact support in $ \Omega $. Since $ \Omega = \cup_{j} \Omega_j $,
\[
C \cup \mbox{Supp} \; u  \subset \Omega_j
\]
for all $ j $ sufficiently large. This implies that the restriction of $ u $ to $ \Omega_j $ belongs to $ H_0(\Omega_j) $ and is admissible for $ C $ for all $ j $ sufficiently large. Also, $ u $ vanishes outside its support and hence $ I(u, \Omega) = I(u, \Omega_j) $ for all $ j $ sufficiently large. As a consequence,
\[
\mbox{Cap}_{\Omega_j}(C) \leq I(u, \Omega_j) = I(u, \Omega)
\]
for all $ j $ sufficiently large. This, in turn, implies that
\[
\ell = \lim_{j \to \infty} \mbox{Cap}_{\Omega_j}(C) \leq I(u, \Omega). 
\]
Since the above inequality holds for every $ u \in H_0(\Omega) $ admissible for $ C $, it follows that $ \ell \leq \mbox{Cap}_{\Omega}(C) $. This establishes the claim.

\medskip

The convergence $ \mbox{Cap}_{\Omega_j}(C) \to \mbox{Cap}_{\Omega}(C) $ ensures that 
\[
\mbox{Cap}_{\Omega_j}(C) < \mbox{Cap}_{\Omega}(C) + \epsilon/2, 
\]
for all $ j $ sufficiently large. Combining the above observation with \eqref{ineq} yields that 
\[ 
\mu_{\Omega_j}(x, x_0) \leq \mbox{Cap}_{\Omega_j}(C) < R 
\]
for all $ j $ sufficiently large. In other words,
$ x \in B_{\Omega_j} (x_0; R) $ for all $ j $ sufficiently large. This establishes \eqref{inclu}.

\medskip

Now, suppose $ \liminf_{j \to \infty} h^K_{\Omega_j}(x_0) < 1/R $ for some $ R > 0 $. It follows from \eqref{inclu} that 
\begin{equation*} 
 B_\Omega(x_0; R - \epsilon) \subset B_{\Omega_j} (x_0; R)   
\end{equation*}
for all $ j $ sufficiently large. Choose $ j $ sufficiently large such that $ h^K_{\Omega_j}(x_0) < 1/R $. Then there is a $ K$-quasiconformal mapping
$ f: \mathbb B^n \to \Omega_j \subset \Omega $ such that $ f(0) = x_0 $
and $ f(\mathbb B^n) \supset B_{\Omega_j}(x_0; R) $. It follows that 
\[
f(\mathbb B^n) \supset B_{\Omega_j}(x_0; R) \supset B_{\Omega}(x_0; R - \epsilon),
\]
and consequently $ h^K_\Omega(x_0)< 1/( R -\epsilon) $. Letting $ \epsilon \ra 0 $, it follows that
$ h^K_\Omega(x_0) \leq 1/R $. Since this inequality holds for every $ R > 0 $ such that $ \liminf_{j \to \infty} h^K_{\Omega_j}(x_0) < 1/R $, letting $ 1/R $ decrease to $ \liminf_{j \to \infty} h^K_{\Omega_j}(x_0) $, it follows that $ \liminf_{j \to \infty} h^K_{\Omega_j}(x_0) \geq h^K_{\Omega}(x_0) $, as required. 
\end{proof}

\begin{thm}
Let $ \Omega \subset \mathbb{R}^n $ be a bounded domain and that $ \{\Omega_j \} $ be a decreasing sequence of domains containing $ \Omega $ such that $\cap_{j} \Omega_j= \Omega $. If $ \mathbb{R}^n \setminus \Omega_1 $ is nonempty and each $ \Omega_j $ satisfies $  \mbox{int}(\overline{\Omega_j})= \Omega_j $, then, for $ K \in [1, \infty) $ fixed, 
\[
h^K_{\Omega}(x) \leq \liminf_{j \to \infty} h^K_{\Omega_j}(x) ,
\]
for every $ x \in \Omega $.  
\end{thm}

\begin{proof}
Fix $ x_0 \in \Omega $ and suppose that $ \liminf_{j\to\infty} h^K_{\Omega_j}(x_0)< 1/r $ for some $ r> 0 $. Passing to a subsequence, we may assume that $ h^K_{\Omega_j}(x_0)< 1/r $ for all $ j $. Then there are $ K$-quasiconformal mappings
$ f_j : \mathbb B^n \to \Omega_j $ with $ f_j(0)=x_0 $
and $ f_j(\mathbb B^n) \supset B_{\Omega_j}(x_0; r) $. Since $ \Omega_j \subset \Omega_1 $ for all $ j $ and $ \mathbb{R}^n \setminus \Omega_1 $ is non-empty, each map $ f_j $ omits two values of $ \overline{\mathbb{R}^n} $ whose spherical distance
is bounded below by a positive constant independent of $ j $. It follows that the family $ \{f_j\} $ is normal. Passing to a subsequence, we may assume that $ \{f_j\} $ converges uniformly on compact subsets of $ \mathbb{B}^n $ to a mapping $ f: \mathbb{B}^n  \to \overline{\Omega}_1 $. 

\medskip

Since $ \Omega \subset \Omega_j $ for every $ j $, it follows that 
\[
B_{\Omega}(x_0; r) \subset B_{\Omega_j}(x_0; r) \subset f_j(\mathbb{B}^n).
\]
As a consequence, 
\[
g_j := f_j^{-1} |_{B_\Omega(x_0; r)} : B_\Omega(x_0; r) \to \mathbb{B}^n 
\]
is a family of $ K$-quasiconformal mappings with $ g_j(x_0) = 0 $. Since $ \mathbb{B}^n $ is bounded, $ \{g_j \} $ is normal and the limit $ g $ is either $ K$-quasiconformal on $ B_\Omega(x_0, r) $ or a constant. Since $ g(x_0) = 0 $, it follows that, in either case, $ g \big( B_\Omega(x_0, r) \big) \subset \mathbb{B}^n $. 

\medskip

It is evident that $ f_j( g_j(x)) = x $ for every $ x \in B_\Omega(x_0; r) $. Since $ g_j(x) \to g(x) \in \mathbb{B}^n $ and $ f_j \to f $ uniformly on compact subsets of $ \mathbb{B}^n $, it follows that $ f \circ g(x) = x $ for every $ x \in B_\Omega(x_0; r) $. As a consequence, $ B_\Omega(x_0, r) \subset f(\mathbb{B}^n) $, which, in turn, implies that $ f $ must be non-constant. In particular, $ f $
is $ K$-quasiconformal on $ \mathbb{B}^n $ and $ f(\mathbb{B}^n) \supset B_\Omega(x_0; r) $.

\medskip

The claim now is $ f(\mathbb{B}^n) \subset \Omega $. To see this, fix $ j_0 \ge 1$ and note that
\[
 f_j(\mathbb{B}^n) \subset \Omega_j \subset \Omega_{j_0} 
 \]
 for all $ j \geq j_0 $. It follows that $ f(\mathbb{B}^n) \subset \overline{\Omega}_{j_0} $. Since $ f(\mathbb{B}^n) $ is an open subset of $ \overline{\mathbb{R}^n} $, it follows that $  f(\mathbb{B}^n) \subset \mbox{int}(\overline{\Omega}_{j_0}) = \Omega_{j_0} $. Since $ j_0 \ge 1$ was arbitrary, it follows that 
 \[
 f(\mathbb{B}^n) \subset \displaystyle\bigcap_{j_0 \ge 1} \; \Omega_{j_0} = \Omega,  
 \]
 and consequently $ h^k_\Omega(x_0) \leq 1/r $.  This yields $ h^K_{\Omega}(x_0) \leq \liminf_{j \to \infty} h^K_{\Omega_j} (x_0) $, as required.
\end{proof}
\subsection{Boundary behaviour of $h^K_{\Om}$} Just as the $K$-squeezing function has a well-defined asymptotic limit at spherically extreme boundary points, so does the $K$-Fridman invariant. The details are similar and, therefore, we will be brief.

\begin{thm}
Let $\Omega \subset \mathbb{R}^n$ be a bounded domain and $\{p_j\}$ be a sequence of points in $\Omega$ converging to $p \in \partial \Omega$. Assume that $ \partial \Omega $ is smooth near $p$ and $p$ is a spherically extreme boundary point. Then,  for $ K \in [1, \infty) $,
\[
h^K_{\Omega}(p_j) \to 0.
\]
\end{thm}

\begin{proof} Proceeding exactly as in the proof of Theorem~\ref{bd-bev-s}, we may assume that $p=0$ and the outward normal to $\partial \Omega$ at the origin points in the direction $(-1,0')$. Following the construction in that proof, we obtain a sequence of domains 
\[
\Omega_j \subset \mathcal{H}=\{x \in \mathbb{R}^n: x_1 > -1\},
\]
 each conformally equivalent to $ \Omega $, such that $0 \in \Omega_j$ for every $j$, $\Omega_j \to \mathcal{H}$ in the local Hausdorff sense. It follows that
 \[
h^1_{\Omega}(p_j)=h^1_{\Omega_j}(0)
\]
for every $ j $. Since $ 0 \leq h^K_{\Omega} \leq h^1_{\Omega} $, it suffices to show that
\[
\limsup_{j \to \infty}h^1_{\Omega_j}(0) \leq 0.
\]
To see this, recall that $\mathcal{H}$ is conformally equivalent to $\mathbb{B}^n$. Let $f: \mathbb{B}^n \to \mathcal{H}$ be a conformal homeomorphism with $f(0)=0$. Since the $ \mu $-metric is conformally invariant and $ \mathcal H $ is conformally equivalent to $ \mathbb B^n $, the relative compactness of $ \mu_{\mathbb B^n}$-balls implies the same property for $ \mu_{\mathcal H} $-balls.
In particular, for any fixed $ R>0 $,
$ B_{\mathcal{H}}(0, R) \subset f(\mathbb{B}^n) = \mathcal{H} $. From this point, an argument similar to the one employed in the first part of Theorem \ref{inc} yields that 
 \[
\limsup_{j \to \infty}h^1_{\Omega_j}(0)<\frac{1}{R}.
\]
Since $R>0$ is arbitrary, it follows that $ \limsup_{j \to \infty}h^1_{\Omega_j}(0)\leq 0
$ as required.
\end{proof}

\section{Relating the $K$-Fridman invariant with the injectivity radius when $n=2$}
Let $ M $ be a $C^1$-smooth Class II Riemannian manifold of dimension $ 2$. By Ferrand's classification, $ M $ is a hyperbolic Riemann surface. Since $ \mu_M $ induces the manifold topology on Class II manifolds, and the Kobayashi metric induces the intrinsic topology of a hyperbolic Riemann surface, these two metrics induce the same topology on $ M $. This observation motivates us to introduce a Kobayashi-metric analogue of the $K$-Fridman function in dimension $2$.

\medskip

Fix $ x \in M $ and let $ \pi: \Delta \to M $ be the universal covering map. Choose $ z_0 \in \pi^{-1}(x) $, and $ A \in \mbox{Aut}(\Delta) $ such that $ A(0)=z_0 $. Since $ \pi \circ A $ is biholomorphic near the origin, for sufficiently small $ r>0 $, the map $ f: \Delta \to M $ defined by $ f(z)=\pi \big(A(rz) \big) $ is $1$-quasiconformal. Moreover, $ f(0)=x $, and $ f(\Delta) $ contains a Kobayashi ball centred at $ x $ since the Kobayashi distance induces the manifold topology. As $ f $ is $1$-quasiconformal, it is $ K$-quasiconformal for every $ K \in [1,\infty) $. 

\medskip

For a hyperbolic Riemann surface $ M $, $ K \in [1, \infty) $ and $ x \in M $, define
\[
H^K_M (x) = \inf_{r \in \mathcal{R}} \frac{1}{r},
\]
where $ \mathcal{R} $ is the set of all $ r > 0 $ for which there exists a $ K$-quasiconformal map $ f : \Delta \to X $ with $ f(0) = x $ and $ f(\Delta) \supset B^k_M(x, r) $, and $ B^k_M(x, r) $ is the Kobayashi ball with center $ x $ and radius $ r $. $H_M^K (\cdot) $ is a biholomorphic invariant since the
Kobayashi metric is itself preserved by such maps. Recall that the Kobayashi metric on the unit disc $ \Delta $ coincides with the hyperbolic distance $ \rho_{\Delta} $.

\medskip

Building on the ideas in the proofs of Theorems $1, 2$ and $3$ of \cite{NTT}, the following results hold. The proofs are exactly the same as in \cite{NTT} and are therefore not repeated.

\begin{prop}
Let $ X $ be a hyperbolic Riemann surface. Then
\[
H_X^K(z) = \inf\{ 1/r : B_X^k(z,r) \mbox{ is simply connected } \},
\]
for all $ z \in X $. 
\end{prop}

\begin{prop}
Let $ D \subset \mathbb{C} $ be a convex domain that contains no complex affine lines and let $ k_D $ denote its Kobayashi distance. Let $ \Gamma $ be a torsion-free discrete group of isometries of $ (D, k_{D}) $, $ X = D/ \Gamma $ and $ \pi: D \to X $ be the quotient map. For $ z \in X $, let $ w $ be any point in $ D $ such that $ \pi(w) = z $. Then
\[
H^K_X (z) \leq \inf_{ \gamma \in \Gamma \setminus \{ \text{Id} \} } \frac{ k_{D} \big( w, \gamma(w) \big)}{2}.
\]
\end{prop}

\begin{prop}
Let $ \Gamma $ be a torsion-free discrete group of isometries of $ (\Delta, \rho_{\Delta}) $, $ X = \Delta/ \Gamma $ and $ \pi: \Delta \to X $ be the quotient map. For $ z \in X $, let $ w $ be any point in $ \Delta $ such that $ \pi(w) = z $. Then
\[
H^K_X (z) = \inf_{ \gamma \in \Gamma \setminus \{ \text{Id} \} } \frac{\rho_{\Delta} \big( w, \gamma(w) \big)}{2}.
\]
\end{prop}

\section{Comparing $s^K$ and $h^K$}

\begin{prop}\label{sq-frid}
Let $ \Omega \subset \mathbb{R}^n $ be a bounded domain, $ n \geq 2 $, and $ K \geq 1 $. Then for every $ p \in \Omega $,
\[
s_\Omega^K (p) \leq \frac{1}{ \gamma_n^{-1} \left( \frac{K}{h_\Omega^K (p)} \right)}
, 
\]
where $ \gamma_n: (1, \infty) \rightarrow (0, \infty) $ is  the $ n$-dimensional Gr\"{o}tzsch modulus function. In particular, if $ \{p_j\} $ is a sequence of points in $ \Omega $ such that
$ s_\Omega^K(p_j) \rightarrow 1, $ then $ h_\Omega^K(p_j)\rightarrow 0 $.

\end{prop}

\begin{proof}
Assume $ s_\Omega^K (p) = r $, and let $ f: \Omega \rightarrow \mathbb{B}^n $ be a $ K$-quasiconformal map that realizes $ s_\Omega^K(p) $, i.e., $ f(p) =0 $ and $ B(0;r) \subset f(\Omega) $. Then $ g: \mathbb{B}^n \rightarrow \Omega $ defined by setting 
$$ g(z) = f^{-1}(rz) $$
is a $ K $-quasiconformal map satisfying $ g(0) = p $. The claim is that the $ \mu_\Omega$-ball
\[
B_\Omega \left( p, \frac{ \gamma_n \big(1/r \big)} {K}\right), 
\]
is contained in $ g(\mathbb{B}^n) $ or equivalently that 
\[
f \left( B_\Omega \Big( p, \frac{ \gamma_n \big(1/r \big)} {K} \Big) \right) \subset B(0;r).
\]
Indeed, if $ q \in B_\Omega \Big( p, \frac{ \gamma_n \big(1/r \big)} {K} \Big) $, then
\[
\gamma_n\left(\frac{1}{ |f(q)|} \right)= \mu_{\mathbb{B}^n} \big( 0, f(q) \big)= \mu_{\mathbb{B}^n} \big( f(p), f(q) \big) \leq K \mu_{\Omega} ( p, q ) < \gamma_n \left(\frac{1}{r} \right).
\]
Since $ \gamma_n $ is strictly decreasing, the above inequality is exactly the assertion that $ |f(q) | < r $. Hence, the claim follows.

\medskip

Next, note that $ g $ is a candidate map in the definition $ h_\Omega^k(p) $ with $ B_\Omega \left( p, \frac{ \gamma_n \big(1/r \big)} {K}\right) \subset g(\mathbb{B}^n) $ and therefore
\[
h_\Omega^K(p) \leq \frac{K} { \gamma_n \big(1/r \big)}.
\]
Since $ r = s_\Omega^K (p)$, it follows that 
\[
\gamma_n \left(\frac{1}{s_\Omega^K (p)} \right) \leq 
\frac{K}{h_\Omega^K (p)}.
\]
The desired inequality then follows by applying the strictly decreasing inverse function $ \gamma_n^{-1} $.

\medskip

Moreover, if $ s_\Omega^K(p_j)\rightarrow 1 $, then evidently $ \frac{1}{s_\Omega^K(p_j)}\rightarrow 1^+ $, and hence 
\[
\gamma_n\left(\frac{1}{s_\Omega^K(p_j)}\right)\rightarrow\infty,
\]
by 
\[
\lim_{t\to1^+}\gamma_n(t)=\infty.
\]
The inequality  
\[
\gamma_n\left(\frac{1}{s_\Omega^K(p_j)}\right)
\leq
\frac{K}{h_\Omega^K(p_j)}.
\]
forces $h_\Omega^K(p_j) \rightarrow 0$.
\end{proof}

The fact that $s^K_{\Om} \ra 1$ implies that $h^K_{\Om} \ra 0$ evidently provides another proof of Theorem $4.8$ on the boundary behaviour of the $K$-Fridman invariant. 

\section{Comments on the Ferrand--Obata theorem}

There are two parts to this concluding section. First, we provide some details of the proof of the domain-version of the Ferrand--Obata theorem, as indicated in $\S 1$. Second, as indicated earlier, the conformal invariants $s^K_{\Om}, h^K_{\Om}$ provide an alternative way to approach this theorem, and we provide the details.

\medskip

For the first part, let $\Om \subset \mbb R^n$, $n \ge 3$, be a $C^2$-smooth bounded domain with a non-compact conformal group ${\rm Conf}(\Om, g_{euc})$. By Liouville's theorem, each
$f \in {\rm Conf}(\Om, g_{euc})$ is a M\"obius map of $S^n$. Since $f$ preserves $\Om$, it must also preserve $\pa \Om$. Since ${\rm Conf}(\Om, g_{euc})$ is non-compact, \cite{GM} shows that every M\"obius group is a convergence group, and therefore, there are points $p, q \in S^n$, with $p= q$ being a possibility, and $f_j \in {\rm Conf}(\Om, g_{euc})$ such that $\{f_j\}$ converges locally uniformly on $S^n \sm \{q\}$ to $p$ and $f_j^{-1}$ converges locally uniformly on $S^n \sm \{p\}$ to $q$.

\medskip

Note that $p, q$ are both in $\pa \Om$. Indeed, pick $a \in \pa \Om \sm \{q\}$. Since $f_j(a) \in \pa \Om$, $j \ge 1$, and $\pa \Om$ is compact, it follows that the limit of $f_j(a)$ which is $p$, must also be in $\pa \Om$. The same argument applied to $f^{-1}_j$ shows that $q$ must also be in $\pa \Om$. It follows that for every $x \in \Om$, $f_j(x) \ra p$.

\medskip

Fix $x_0 \in \Om$, write $p_j = f_j(x_0)$ and note that $p_j \ra p \in \pa \Om$. Rescale the boundary $\pa \Om$ near $p$ by the affine conformal maps $A_j, T_j$, as explained in Theorem \ref{bd-bev-s}, to see that the scaled domains $\Om_j = T_j \circ A_j \circ f_j(\Om)$ converge locally uniformly in the Hausdorff metric to the half space $\mathcal H$ as before. Furthermore, $\{T_j \circ A_j \circ f_j\}$ is a family of M\"obius maps that satisfy $T_j \circ A_j \circ f_j(x_0) = 0 \in \mathcal H$ for $j \ge 1$ and the $\Om_j$'s do not intersect a uniform open set in $\mbb R^n$ due to the fact that the boundaries $\pa \Om_j$ are locally close to $\pa \mathcal H$ for $j$ large. Thus, $\{T_j \circ A_j \circ f_j\}$ is normal. 

\medskip

Again, by \cite{GM}, there are two possibilities. First, there are two points $\al, \be \in S^n$, both possibly the same, such that $T_j \circ A_j \circ f_j$ converges locally uniformly on $S^n \sm \{\be\}$ to $\al$ and the inverses $(T_j \circ A_j \circ f_j)^{-1}$ converge locally uniformly on $S^n \sm \{ \al\}$ to $\be$. Second, $T_j \circ A_j \circ f_j$ converges locally uniformly on $S^n$ to a M\"obius limit. Using similar arguments as above on the relative location of $\al, \be$, it can be seen that the first possibility cannot arise. The limit of $T_j \circ A_j \circ f_j$ is then a M\"obius equivalence between $\Om$ and the half space $\mathcal H$, which, in turn, is conformally equivalent to $\mbb B^n$.

\medskip

For the second part, here are two variants of the domain-version of the Ferrand--Obata theorem whose proofs rely on the conformal invariants $s^K_{\Om}, h^K_{\Om}$. Their efficacy derives from the fact that we do not have to appeal to either Liouville's theorem or other deeper properties of M\"obius groups. The mere presence of a divergent sequence of conformal automorphisms that converges to a spherically extreme smooth boundary point is sufficient to deduce conformal equivalence with $\mbb B^n$.

\begin{thm}\label{WR-type}
Let $\Omega \subset \mathbb{R}^n$ be a bounded domain. Suppose that there exists an $x_0 \in \Omega$ and a sequence $\{f_j\}$ of conformal automorphisms of $\Omega$ such that $f_j(x_0) \to p \in \partial \Omega$. If $\partial \Omega$ is smooth near $p$ and $p$ is spherically extreme, then $\Omega$ is conformally equivalent to the unit ball $\mathbb{B}^n$.
\end{thm}
\begin{proof}
Let $p_j=f_j(x_0)$. Then by Theorem~\ref{bd-bev-s},
\[
s^1_{\Omega}(x_0)=s^1_{\Omega}(p_j) \to 1
\]
It follows from Proposition~\ref{s-ball} that $\Omega$ is conformally equivalent to $\mathbb{B}^n$. The same conclusion is obtained by using the boundary behaviour of $h^K_{\Om}$.
\end{proof}

\begin{thm}
Let $\Omega \subset \mathbb{R}^n$ be a smoothly bounded domain such that the group of conformal automorphisms of $\Omega$ acts transitively on $\Omega$. Then $\Omega$ is conformally equivalent to $\mathbb{B}^n$.
\end{thm}
\begin{proof}
Note that the farthest point on $\pa \Om$ from the origin, say $p$, is a spherically extreme boundary point. Let $x_0 \in \Omega$. By hypothesis, there exists a sequence $\{f_j\}$ of conformal automorphisms of $\Omega$ such that $f_j(x_0) \to p$. The proof now follows from Theorem~\ref{WR-type}.
\end{proof}

Notice that the last two theorems are in the spirit of the Wong--Rosay theorem in complex analysis of
several variables.

\section{Questions}

First, as seen in $\S 3.5$, there is a critical dilatation $K_0 > 1$ at which the family of invariants $\{s^K_{\mathcal E}: K \ge 1 \}$ exhibits a transition. For each $1 \le K < K_0$, $s^K_{\mathcal E}$ is non-constant and approaches $1$ at the boundary $\pa \mathcal E$. It would be interesting to calculate $s^K_{\mathcal E}$ for $K < K_0$. Note that the conformal group of an ellipsoid is compact unless it is a ball.

\medskip

Second, the boundary behaviour of $s^K_{\Om}, h^K_{\Om}$ has been studied near spherically extreme smooth boundary points. We believe that the assumption of spherically extreme is not needed and that smoothness alone should suffice to obtain that $s^K_{\Om} \ra 1$ and $h^K_{\Om} \ra 0$, but we have not been able to show this. The rationale for this expectation is that, in the quasiconformal category, smoothness of the boundary already forces the local model near a boundary point to be a half-space. This does not happen in the holomorphic setting where Levi geometry determines the local model domain.

\medskip

Finally, the invariants $s^K_{\Om}, h^K_{\Om}$ encode a metric in terms of which the largest ball is measured; it is the euclidean metric for $s^K_{\Om}$, while for $h^K_{\Om}$, it is the $\mu$-metric. As seen in $\S 5$, the $K$-Fridman invariant can also be defined in terms of the Kobayashi metric, and this allows it to be connected to the injectivity radius. It might be worthwhile to define and study these invariants in terms of other metrics to connect them with other geometric objects. This could have other applications.

\section{Concluding Remarks}

The holomorphic Squeezing function and Fridman invariant, which are well known biholomorphic invariants, are shown to have suitable analogues in the quasiconformal setting. The presence of a dilatation factor $K \ge 1$ allows the construction of an ordered heierarchy of such invariants, which adds to their utility. Although this family of invariants shares some properties, examples show that they retain their peculiarities in general. Despite this, they are related in a computable manner. One of their strengths derives from the fact that they provide a common approach to explore domain-variants of the Ferrand--Obata and Wong--Rosay theorems.

\medskip

It is likely that the $K$-Squeezing function and $K$-Fridman invariants will prove to be useful in other contexts and some questions that arose in this investigation have been outlined above for the interested reader. Surely, there are others. We hope to return to such questions later.

\medskip
\noindent {\bf Declaration}: No AI models were used to generate, verify, or draft any of the mathematical content or text presented in this work. 


\end{document}